\documentclass[12pt,eqno,sumlimits]{amsart}

\usepackage{amssymb, amscd, amsmath, epsfig, mathtools}
\usepackage[utf8]{inputenc}
\usepackage{amsthm}
\usepackage{enumerate}
\usepackage{xcolor}
\usepackage{scalerel}
\usepackage{soul}
\usepackage{tikz-cd}
\usepackage{esint}
\usepackage{slashed}
\usepackage[textsize=tiny]{todonotes}
\usepackage{pinlabel} % for labelling the figure
\usepackage{mathtools}
\usepackage{url}
\usepackage{hyperref}
\usepackage{cancel}

\newtheorem{innercustomgeneric}{\customgenericname}
\providecommand{\customgenericname}{}
\newcommand{\newcustomtheorem}[2]{%
  \newenvironment{#1}[1]
  {%
   \renewcommand\customgenericname{#2}%
   \renewcommand\theinnercustomgeneric{##1}%
   \innercustomgeneric
  }
  {\endinnercustomgeneric}
}

\newcustomtheorem{customthm}{Theorem}
\newcustomtheorem{customlem}{Lemma}
\newcustomtheorem{customprop}{Proposition}
\newcustomtheorem{customcor}{Corollary}

\DeclareMathSymbol{\hy}{\mathord}{AMSa}{"39}

\usepackage[margin=1.0in]{geometry}

\newtheorem{thm}{Theorem}[section]
\newtheorem*{theorem*}{Theorem}

\newtheorem{cor}{Corollary}[section]
\newtheorem{lem}{Lemma}[section]
\newtheorem{prop}{Proposition}[section]
\theoremstyle{definition}
\newtheorem{rem}{Remark}[section]
\newtheorem{exm}{Example}[section]

\newcommand{\C}{\mathbb C}
\newcommand{\R}{\mathbb R}

\newcommand{\Z}{\mathbb Z}

\newcommand{\tr}{{\rm tr}}
\newcommand{\rr}{{\rm R}}
\newcommand{\rk}{{\rm K}}
\newcommand{\re}{{\rm RE}}
\newcommand{\RE}{{\rm RE}}
\newcommand{\rb}{\text{B}}
\newcommand{\dg}{\dot \gamma}

\usepackage{mathtools}

\newcommand{\sr}{{\rm S}}
\newcommand{\crr}{{\rm CR}}

\newcommand{\diff}{{\rm Diff}}
\newcommand{\ism}{{\rm Isom}}

\newcommand{\PP}{\mathbb P}
\newcommand{\sgn}{{\rm sgn}}
\newcommand{\modd}{{\rm mod}}

\begin{document}

\title
{The geometry of loop spaces IV:  Closed Sasakian manifolds} 
\author[Y. Maeda]{Yoshiaki Maeda}
\address{Tohoku Forum for Creativity, Tohoku University}
\email{yoshimaeda@tohoku.ac.jp}
\author[S. Rosenberg]{Steven Rosenberg}
\address{Department of Mathematics and Statistics, Boston University}
\email{sr@math.bu.edu}

\begin{abstract}
We prove that a closed regular $(4k+1)$-Sasakian manifold $(M,h_0)$ admits a family of metrics $h_\rho, \rho\geq 0,$  such that 
$\pi_1({\rm Isom}(M, h_\rho))$, the fundamental group of the isometry group, is infinite for $\rho>0.$  For $M= S^{4k+1}$, this result holds for all $\rho>0$, but fails at $\rho=0.$
\end{abstract}

\maketitle

\section{Introduction}
It is well known that the dimension of the isometry group $\ism(M^n,g)$ of a closed Riemannian manifold $(M^n,g)$ is at most $n(n+1)/2,$ with the maximum achieved for the standard metric on $S^{n}.$  Nevertheless,  this large isometry group has small fundamental group: 
$\pi_1(\ism(S^{n+1}, g_{\rm std})) = \Z_2.$  (Throughout, all fundamental groups are based at the identity map.) In contrast, for 
a generic  $(M,g)$, we expect that  $\ism(M,g)$ is trivial. Thus it seems interesting to produce examples, other than flat tori,  where $\pi_1(\ism(M,g))$ is infinite.  

In this paper, we construct 
a one-parameter family of distinct metrics $h_{\rho}$ on any closed Sasakian manifold  $M^{4k+1}$ such that $\pi_1(\ism(M, h_\rho))$ is infinite.
\begin{customthm}{5.1}\label{thm:5.1}
 Let 
 $(M, g, \phi , \xi , \eta )$ be a $(4k+1)$-dimensional closed
 regular  Sasakian manifold, where $k \ge 1$. 
The metric
 $h_\rho = g + \rho^2 \eta \otimes \eta$
has $|\pi_1(\ism(M, h_\rho)| = \infty$
 for $\rho\in \R^+, \rho \gg 0.$  For $\rho \neq \rho',$ $h_\rho$ and $h_{\rho'}$ are not isometric.
\end{customthm}

Theorem \ref{thm:5.1} will be proved in \S5.  We note that the proof does not work for 
dim$(M) = 4k+3.$

For the case of $S^{4k+1}$, 
we refine Theorem \ref{thm:5.1}:

\begin{customthm}{6.1}\label{thm:6.1}
Let $k \ge 1$. 
Let $ g_{\rm std}$  be the standard metric on $S^{4k+1}$,
with $\eta$  the contact form on
$S^{4k+1}$ associated with the Hopf fibration. 
For the metric $h_\rho = g_{\rm std} +\rho^2 \eta \otimes \eta $,
$|\pi_1(\ism(S^{4k+1},h_\rho))| = \infty$ for 
$\rho^2 \neq 0$.  
\end{customthm}

Thus the standard metric with small $\pi_1$ is the endpoint of a family with infinite $\pi_1.$

Using these techniques, in \S7 we obtain a similar result for an infinite dimensional group.

\begin{customthm}{7.1}\label{thm:7.1a}
Let $\diff_{\eta,str}(S^{4k+1})$ be the group of strict contact diffeomorphisms  of $S^{4k+1}.$
Then
$|\pi_1(\diff_{\eta,str} (S^{4k+1}))| = \infty$.
\end{customthm}

In fact, a  stronger result has been given in \cite{casals}; we reproduce this result in \cite{MRT6}.

As an outline of the paper, in \S2, we recall the setup of Sasakian manifolds.  In \S3, we outline the proof of Thm.~\ref{thm:5.1}, which is based on the Wodzicki-Chern-Simons classes on the loop space $LM$ used in
\cite{MRT5,MRT4}. 
In \S4, we introduce $h_\rho$ and compute its curvature. In \S5 we prove that the pullback of the Wodzicki-Chern-Simons class on  $LM$ to a top class on $M$ is nonvanishing for $\rho\gg 0.$  This involves detailed calculations. 
By the results of \S3, we obtain Thm.~\ref{thm:5.1}. In \S6, we prove Thm.~\ref{thm:6.1}, using the special nature of the standard metric on the sphere.  We obtain similar results for other Sasakian space forms.  
In \S7, we prove Thm.~\ref{thm:7.1}.  Technical results are proven in Appendices.

\section{Sasakian manifolds}
We first  recall the definition of a Sasakian manifold. We work in local coordinates in preparation for later calculations; coordinate-free versions of the definitions are in \cite{Blair}.

A $(2k+1)$-dimensional manifold $M$ is called
an {\it almost contact manifold} if there is an almost
contact structure $(\phi, \xi, \eta ) $ consisting
of a $(1,1)$-tensor field $\phi$, a {\it characteristic vector field}
$\xi$, and a $1$-form $\eta$ satisfying
\begin{equation}\label{eq:one}
 (\phi^2)_i{}^j = \phi_i{}^k \phi_k{}^j
        = -\delta_i{}^j + \eta_i \xi^j ,
\quad \eta_k \xi^k =1 , 
\quad \phi_{j}{}^{k} \xi^j =0,
\quad \phi_{j}{}^k \eta_k =0. 
\end{equation}
($(M, \eta)$ is a {\it contact manifold} if $\eta\wedge (d\eta)^k\neq 0$ pointwise.)
Let $g$ be a Riemannian metric on 
$M$ associated with the almost contact structure
$(\phi, \xi, \eta)$, {\it i.e.,} 
\begin{equation}\label{eq:two}
g_{ij} = g_{k\ell} \phi_i{}^k \phi_j{}^\ell 
   + \eta_i \eta_j.    
\end{equation}
Such a metric always exists. Then $(g,\phi,\xi, \eta) $ 
is called {\it an almost contact metric structure} 
for the almost contact structure
$(\phi, \xi , \eta)$. 
We set 
\begin{equation}
\label{associative-metric}
\phi_{ij} = \phi_i{}^kg_{kj},
\qquad 
\phi^h{}_j = g^{h \ell}
\phi_\ell{}^kg_{kj} .  
\end{equation}
We note that (\ref{eq:two}) implies 
\begin{equation}\label{eq:dual}
g_{ij}\xi^j = \eta_i,
\end{equation}
so  $\xi$ and $\eta$ are dual.
Similarly, $\phi$ is skew-symmetric, since
$$\phi_{ji}  = \phi_j{}^s g_{si} = 
\phi_{j}{}^s (g_{k\ell}\phi_s{}^k 
\phi_i{}^\ell + \eta_s \eta_i) 
=g_{k\ell}(- \delta_j{}^k +\eta_j \xi^k)
\phi_i{}^k
=-g_{\ell}\phi_i{}^{\ell} + \eta_j \eta_{\ell} \phi_i{}^\ell
=-\phi_{ij}. $$

An almost contact metric structure is a 
{\it contact metric structure} (and $(M,g,\phi, \xi, \eta)$ is 
a {\it contact metric manifold}) if 
$\phi =- d\eta$, 
or equivalently, 
\begin{equation}\label{eq:con}\phi_{ij} = -\frac{1}{2}( \nabla_i \eta_j -\nabla_j \eta_i)
 := -\frac{1}{2}( (\nabla_i \eta)_j -(\nabla_j \eta)_i),
\end{equation}
where $\nabla$ is the Levi-Civita connection  for $g$. 
Finally, a 
contact metric manifold 
$(M, g, \phi, \xi, \eta)$ 
is {\it Sasakian} (or {\it normal}) if 
\begin{equation}\label{(A.2)}
    \nabla_i \phi_j{}^k 
                = g_{ij}\xi^k - \eta_j \delta_i{}^k.
\end{equation} 

As in \cite[p.~87]{Blair}, writing (\ref{(A.2)}) as $(\nabla_X\phi)(Y) = g(X,Y)\xi -\eta(Y)(X)$ and setting $Y=\xi$ gives
$-\phi(\nabla_X\xi) = (\nabla_X\phi)(\xi) = \eta(X)\xi -X.$  Apply $\phi$ to both sides, where 
$\phi^2(X) = -X + \eta(X)\xi$, and use
$0 = g(\nabla_X\xi, \xi) = \eta(\nabla_X\xi)$ to obtain $\nabla_X\xi = -\phi(X).$  In coodinates, 
this becomes  $\nabla_i\xi^j = -\phi_i^j$, or equivalently
\begin{equation}\label{eq:refine}
-\phi_{ij} = \nabla_i\eta_j.
\end{equation}
This refines (\ref{eq:con}) for Sasakian manifolds.

Since $\phi_{ij}$ is skew-symmetric, (\ref{eq:refine}) implies 
\begin{equation}\label{eq:ek}\nabla_i \eta_j + \nabla_j \eta_i =0.
\end{equation}
Therefore, $\xi^i =g^{ik} \eta_k$  is a unit Killing vector field on
$M$ which generates a one-parameter group of isometries.

The following result is standard.
\begin{lem}\label{lem:one}
If $(M,g, \phi, \xi, \eta)$ is a $(2k+1)$-dimensional contact metric manifold,
then
\begin{equation}
 \eta \wedge (d\eta)^{k} = (-1)^k\eta\wedge\phi^k \neq 0.    
\end{equation}
$\eta \wedge (d\eta)^{k}$ is the volume form on $M$, so $\int_M \eta \wedge (d\eta)^{k}\neq 0.$ Thus $M$ is a contact manifold.
\end{lem}
\begin{proof}
By (\ref{eq:one}) and (\ref{eq:two}), 
for 
$m \in M$ 
there is an orthonormal basis $\{e_1,\cdots, e_{2k+1}\}$
on $T_mM$ 
such that 
$e_1= \xi, \,  \phi e_{1+\alpha}= e_{k+1+\alpha},$
$     \phi e_{k+1+\alpha}= -e_{1+ \alpha},$ for 
$\alpha= 1, \cdots , k.$
Thus, 
\begin{align*}
 \eta\wedge (d\eta)^k 
    (e_1, e_2\cdots, e_{k+1}, e_{k+2}, \cdots, e_{2k+1} )
&    =  \eta (e_1) g(\phi e_2, e_{k+1}) \cdots 
        g(\phi e_{k+2}, e_{2k+1}) = 1 \neq 0. 
\end{align*}
\end{proof}

A basic example of a Sasakian manifold is 
as follows. 
\begin{exm}\label{exm:1}
If $M=S^{2k+1}$, the $(2k+1)$-dimensional standard sphere of radius $1$, then
$M$ is a Sasakian manifold with respect to 
the Hopf fibration 
\begin{equation}
\label{(4)}
\pi : S^{2k+1} \to \C\PP^k,    
\end{equation}
where the characteristic vector field
 $\xi$ is the unit vector along the fiber, and  $\eta(\xi) = 1$ by (\ref{eq:one}).  
\end{exm}
\medskip

In this paper, we assume  that $M$ is closed and {\em regular}, 
{\it i.e.,}
the integral curves characteristic vector field $\xi$ starting at any point $m\in M$  are periodic.
 By \cite[p.~37]{Blair}, the period of the integral curves is independent of $m$.  In this case, we have the circle fibration $S^1\to M \to M/S^1$, and $M/S^1$ is a $2k$-dimensional symplectic manifold.   
 The example of the Hopf fibration is a special case.  If $\xi$ is regular, then there is a circle action
 \begin{equation}\label{action}a:S^1\times M\to M,\ a(\theta, m) = R_\theta(m),
 \end{equation}
 given by rotation by angle $\theta$ in the circle fiber. This action  is via isometries for the metrics $h_\rho$  below (see Cor.~\ref{cor:hrho}).

The Theorems in the Introduction entail proving that the loop of isometries
$a^I:S^1\to\ism(M,h_\rho)$, $a^I(\theta)(m) = a(\theta,m)$, has nontrivial class
$[a^I]\in \pi_1(\ism(M, h_\rho))$ and that the powers $[a^I]^n$ are all distinct.

\section{Outline of the proofs}

 As explained in \cite[\S2.2]{MRT4}, a Riemannian metric $g$ on $M$ determines the Wodzicki-Chern-Simons form $CS^W_g\in \Lambda^{2k-1}(LM)$ given by 
\begin{align}\label{eq:CSW}\MoveEqLeft{}{CS^W_{2k-1}(g)(X_1,...,X_{2k-1}) }\\
&= 
\frac{k}{2^{k-2}} 
\sum_{\sigma} {\rm sgn}(\sigma) \int_{S^1}\tr[
 (R(X_{\sigma(1)},\cdot)\dg)
 (\Omega^M)^{k-1}(X_{\sigma(2)},\ldots,X_{\sigma(2k-1)} )],\notag
 \end{align}
at the loop $\gamma\in LM$. Here $\sigma\in (1,\ldots,2k-1)$, the permutation group on 
$\{1,\ldots,2k-1\},$ $R$ is the curvature tensor of $g$,  $\Omega$ is the curvature two-form of $g$, and the $X_i\in T_\gamma LM$ are vector fields $X_{i,\theta}\in T_{\gamma(\theta)}M$ defined along $\gamma$.  This form comes from replacing the usual trace in the definition of Chern-Simons classes with the Wodzicki residue, which is a trace on the Lie algebra of the structure group for natural metric connections on $TLM.$

Associated to the action $a$ is the map $a^L:M\to LM$ from $M$ to the free loop space $LM = C^\infty(S^1, M)$ given by $a^L(m)(\theta) = a(\theta,m).$ This produces a top form $a^{L,*}CS^W \in \Lambda^{2k-1}(M).$
In \cite[Appendix A]{EMS}, we prove 
\begin{prop}\label{prop:app}
\begin{equation}\label{eq:2.1}\int_M a^{L,*}CS^W_{2k-1}(g) \neq 0 \Rightarrow 0\neq [a^I]\in \pi_1(\ism(M,g))
\ {\rm and}\ [a^I]^n \neq [a^I]^m\ {\rm for} \ m\neq n.\,
\end{equation}
for $g = h_\rho.$  
\end{prop}
The proof explicitly uses that the action is via isometries.

By the long  calculation in \S5 and Appendix \ref{app:prop1}, $a^{L,*}CS^W_{2k-1}(h_\rho)$ is a form-valued polynomial in $\rho$
whose leading coefficient is a nonzero multiple of $\eta\wedge (d\eta)^k$ (see (\ref{eq:a*CSW}) and Prop.~\ref{prop:1}). 
Thus   $\int_M a^{L,*}CS^W_{2k-1}(h_\rho)$ is a polynomial in $\rho$
whose leading coefficient is a nonzero multiple of the $h_\rho$-volume of $M$. Therefore,
$\int_M a^{L,*}CS^W_{2k-1}(h_\rho)\neq 0$ for $\rho \gg 0,$ so  Thm.~\ref{thm:5.1} follows from (\ref{eq:2.1}).  For the standard metric on $S^{4k+1}$, the form-valued polynomial consists of only one term, so $\int_M a^{L,*}CS^W_{2k-1}(h_\rho)\neq 0$ for all $\rho\neq 0$, which implies 
Thm.~\ref{thm:5.1}.

The tensorial calculation in \S5 analyzes the integrand in (\ref{eq:CSW}), using the curvature calculations for $h_\rho$ in \S4.  The integrand initially consists of an overwhelming number of terms, but the Bianchi identity and other symmetries are used repeatedly to eliminate most terms.

\section{Curvature formulas for the metric $h_\rho$}
Let $(M,g,\phi, \eta, \xi)$ be a regular Sasakian manifold. Define 
 the new metric $h = h_\rho$ by 
\begin{equation*}
h_{ij} = g_{ij}   +\rho^2 \eta_i \eta_j.
\end{equation*}
To motivate this metric, we recall from \cite[\S3.2]{MRT4} that a regular contact manifold $M$ is the total space of an $S^1$-bundle $L_k\stackrel{\pi}{\to} N$ over a Riemannian symplectic manifold $(N, \omega, g_N)$ with $[\omega]\in H^2(N,\Z)$ and $c_1(L) = k\omega$ for our choice of $k\in \Z.$  $L_k$ has the contact metric
$g_k = \pi^*g_N + \eta_k\otimes \eta_k$, where $\eta_k$ is dual to $\xi_k$, the unit tangent vector field to the fibers.  $L_1\stackrel{\alpha}{\to} L_k$ is a $k$-fold cover, so $g_k$ pulls back to a locally isometric metric on $L_1.$

Noting that $\eta_1 = k\alpha^*\eta_k$, we have a $\Z$-family of metrics on $L_1$ of the form $\alpha^*g_k + k^{2}\eta_1\otimes \eta_1$.  This is a special case of $h_\rho.$

We will compute the Levi-Civita connection and the Riemannian curvature for $h.$  Some proofs are in Appendix \ref{app:curv}.

\begin{lem}
Let $h^{ij}$ be the inverse of the metric $h_{ij}$.
Then 
\begin{equation*}
 h^{ij}= g^{ij} + \alpha \xi^i \xi^j\ {\rm for}\ \alpha = -\frac{\rho^2}{1+\rho^2}.   
\end{equation*}
\end{lem}

\begin{proof}

\begin{align*}
 h_{ik} h^{kj}&= (g_{ik} + \rho^2 \eta_i \eta_k)
                  (g^{kj} + \alpha \xi^k \xi^j) 
 = \delta_i^j +(\rho^2 + \alpha + \rho^2 \alpha) \eta_i \xi^j 
 = \delta_i^j,
\end{align*}
for $\alpha= - \frac{\rho^2}{1+\rho^2}$. 
\end{proof}

Let $\Gamma_{ij}^{k}$ and $\bar{\Gamma}_{ij}^{k}$
be the Christoffell symbols of $g$ and
$h$ respectively:   
\begin{align*}
\Gamma_{ij}^{k}  
&= 
\frac{1}{2} g^{k\ell}
     (\partial_i g_{j\ell} + \partial_j g_{i\ell} 
       -\partial_\ell g_{ij} ),\ 
 \bar{\Gamma}_{ij}^{k} 
= \frac{1}{2} h^{k\ell}
     (\partial_i h_{\ell j} + \partial_j h_{i\ell} 
       -\partial_\ell h_{ij} ). 
       \notag
\end{align*}

\begin{lem}\label{lem:3.2}
\begin{equation*}
  \bar{\Gamma}_{ij}^{k}  
    =  \Gamma_{ij}^{k}   
     - \rho^2 (\phi_i^k \eta_j + \phi_j^k \eta_i). 
\end{equation*}    
\end{lem}

\begin{proof}  See Appendix \ref{app:curv}.
\end{proof}

We want the action (\ref{action}) to be an action by isometries for $h_\rho.$ Since $\Vert \xi\Vert^2_{h_\rho} = 1+ \rho^2$ by (\ref{eq:one}), (\ref{eq:two}) and the definition of $h_\rho,$  $\bar\xi := (1+\rho^2)^{-1/2}\xi$ has $\Vert \bar\xi\Vert^2_{h_\rho} = 1.$  We can consider (\ref{action}) to be an action with $(d/d\theta) a(\theta,m) = \bar\xi$ by modifying the $g$-length of the orbits. If 
orbits have $g$-length $\lambda$, then these orbits now have $h_\rho$-length also equal to $\lambda.$

\begin{cor}\label{cor:hrho}  The $S^1$-action $\{R_\theta\}$ in (\ref{action}) is an action via isometries for $h_\rho.$
\end{cor}

\begin{proof}
       We have
   $\bar\eta = \bar\eta_i dx^i$ is given by
   $$\bar\eta_i = (h_\rho)_{ij}\bar\xi^j = (g_{ij} + \rho^2\eta_i\eta_j)((1+\rho^2)^{-1/2} \xi^j)
   = (1+\rho^2)^{1/2}\eta_i
  .$$
Let $\bar\nabla$ be the Levi-Civita connection for $h_\rho.$ By Lem.~\ref{lem:3.2}, 
$$\bar\nabla_i\bar\eta_j := (\bar\nabla \bar \eta)_i = 
(1+\rho^2)^{1/2}[\nabla_i\eta_j + 
   \rho^2(\phi_i{}^s \eta_j + \phi_j{}^s \eta_i)\eta_s]
=(1+\rho^2)^{1/2}(\nabla_i\eta_j),$$
so 
$$\bar\nabla_i\bar\eta_j +  \bar\nabla_j\bar\eta_i 
  =(1+\rho^2)^{1/2}(\nabla_i\eta_j + \nabla_j\eta_i)=0.$$

Therefore, $\bar\xi$ is an $h_\rho$-unit Killing vector field and the rotation action is via isometry.
\end{proof}

We now compute the curvature tensor for the metric $ h_\rho$.

\begin{lem}\label{lem:3.3}
Let $\bar{R}_{kji}{}^h$ and $R_{kji}{}^h$
be the curvature tensor of the
Riemannian metrics $h_\rho$ and $g$ respectively.  
Then
\begin{align}\label{eq:14a}
\bar{R}_{kji}{}^h &= R_{kji}{}^h 
 - \rho^2 (\phi_{ki} \phi_j{}^h -\phi_k{}^h \phi_{ji}
+  2 \phi_{kj} \phi_i{}^h 
+2\eta_k \eta_i \delta_j{}^h -2\eta_j \eta_i \delta_k{}^h 
+g_{ki} \eta_j \xi^h -g_{ji}\eta_k \xi^h ) \notag
\\
&\qquad - \rho^4(\eta_k \eta_i\delta_j{}^h  
- \eta_j \eta_i\delta_k{}^h ). 
\end{align}
\end{lem}

\begin{proof}  This is in Appendix \ref{app:curv}. 
\end{proof}

\section{Computation of $CS^W$}

In this section, we compute the Wodzicki-Chern-Simons form for $(M, h_\rho),$  or more precisely its pullback $(a^L)^*CS^W$ (\ref{eq:a*CSW}). As explained in Prop.~\ref{prop:app}, the nonvanishing of 
$(a^L)^*CS^W$ (Cor.~\ref{cor:a4k}) is the key step in the proofs of 
 the main Thms.~\ref{mainthm}, \ref{thm:sphere}.
\bigskip

\noindent {\bf Notation:} We   
introduce some local coordinates notation.  
Let ${\mathfrak{T}}_P^Q (M)$ 
and $\Lambda^P(M)$
be the space  of tensor fields of type 
$(Q,P)$ and $P$-froms on $M$ respectively. 
For an $(0,P)$-tensor 
$S_{i_1 \ldots i_P}$ and  $(0,P^\prime)$-tensor 
$S^\prime_{i_1 \ldots i_P}$, its tensor product is
\begin{equation*}
(S\otimes S^\prime)_{i_1 \cdots i_{P+P^\prime}} 
:= S_{i_1 \cdots i_{P}} \cdot S^\prime_{i_P \cdots i_{P+P^\prime}}
\in \mathfrak{T}^0_{P+P'}(M).
\end{equation*}
For an $(0,P)$-tensor $S_{i_1 \ldots i_P}$, 
its skew-symmetrization with respect to
the indices $i_1, \ldots, i_P$ is denoted by 
\begin{equation*}
\label{skew-symmetric}
S_{[i_1 \cdots i_P]}
:= \frac{1}{P!} 
\sum_{\sigma \in (1, 2, \ldots, P)}
\sgn(\sigma) 
S_{i_{\sigma(1)} \ldots i_{\sigma(P)}} \in \Lambda^{P}(M).
\end{equation*}   
For $A \in \Lambda^P (M)$ and 
$A^\prime \in \Lambda^{P^\prime} (M)$,
 the exterior product is 
\begin{align}
\label{wedge-product}
(A\wedge A^\prime)_{i_1 \cdots i_{P+P^\prime}} 
&= A_{[i_1 \cdots i_{P}} 
   A^\prime_{i_{P+1} \cdots i_{P+P^\prime}]}
   \\
&:=\frac{1}{(P+P^\prime)!}
\sum_{\sigma \in (1, \cdots, P+P^\prime)}
\sgn(\sigma) 
A_{i_{\sigma (1)} \cdots i_{\sigma (P)}} 
   A^\prime_{i_{\sigma(P+1)} \cdots 
     i_{\sigma(P+P^\prime) }}.
     \notag
\end{align}
 $\tr$ means the contraction of the superscript index and the subscript index of the $\mathfrak T_1^1$ part of a tensor.
For example, for
$S=S_{i_1 \cdots i_P k}{}^\ell \in {\mathfrak{T}}_P^0 (M) 
     \otimes {\mathfrak{T}}_1^1 (M)$, 
\begin{equation}\label{trace}
\tr( S)_{i_1 \cdots i_{P} }   
= S_{i_1 \cdots i_P k}{}^k.
\end{equation}
This notation extends to 
a ``tensor and contract" operation on ${\mathfrak{T}_P^0} \otimes {\mathfrak{T}}_1^1$.
Namely, for 
$S_{i_1 \cdots i_P}{}_{\ell^\prime} {}^{\ell_\alpha} \in {\mathfrak{T}_P^0} \otimes {\mathfrak{T}}_1^1$
and 
$S^\prime_{i_1 \cdots i_{P^\prime}}{}_{\ell_\beta}{}^{\ell^\prime} 
      \in {\mathfrak{T}_{P'}^0} \otimes {\mathfrak{T}}_1^1$,
set
\begin{equation*}
\
(S \otimes S^\prime)_{i_1 \cdots i_{P+P^\prime}}{}_{ \ell_\beta}^{\ \ \ell_\alpha}   
= S_{i_1 \cdots i_P}{}_{\ell^\prime}{}^{\ell_\alpha}
S^\prime_{i_1 \cdots i_{P^\prime}}{}_{\ell_\beta}{}^{\ell^\prime}\in \mathfrak{T}_{P+P'}\otimes\mathfrak{T}^1_1.
\end{equation*}
This similarly extends to a ``wedge and contract" operation on 
$\Lambda^P (M) \otimes {\mathfrak{T}}_1^1$:
For 
$A=A_{i_1 \cdots i_P \ell'}{}^{\ell_\alpha}
\in \Lambda^P(M) 
     \otimes {\mathfrak{T}}_1^1 (M)$, 
and 
$B_{i_1 \cdots i_{P'} \ell_\beta}{}^{\ell'} \in \Lambda^{P^\prime} (M) 
     \otimes {\mathfrak{T}}_1^1 (M)$,  set 
\begin{align*}
\MoveEqLeft{(A\wedge B)_{i_1 \cdots i_{P+P^\prime}\ell_\beta}{}^{  \ell_\alpha} 
 = (A\otimes B)_{[i_1 \cdots i_{P+P^\prime}]\ell_\beta}{}^{  \ell_\alpha} }\\ 
 &= 
 \frac{1}{(P+P^\prime)!}
\sum_{\sigma \in (1, \cdots, P+P^\prime) }
\sgn(\sigma) 
A_{i_{\sigma (1)} \cdots i_\sigma(P)}{}_{\ell^\prime}{}^{\ell_\alpha}
B_{i_{\sigma (P+1)} \cdots 
    i_{\sigma (P+P^\prime})}{}_{\ell_\beta}{}^{\ell^\prime}\\
    &\qquad \in \Lambda^{P+P'}(M) \otimes \mathfrak{T}_1^1.
\end{align*}
For 
$S_{i_1 \cdots i_P} , {S^\prime}_{i_1 \cdots i_P}
  \in {\mathfrak{T}}_P^0$,
define an equivalence relation by  
\begin{equation*}
S_{i_1 \cdots i_P} \sim {S^\prime}_{i_1 \cdots i_P}  
\ \modd (i_1 \cdots i_{P} )
\Leftrightarrow
S_{[i_1 \cdots i_P]} = {S^\prime}_{[i_1 \cdots i_P]},     
\end{equation*}
and for 
$S_{i_1 \cdots i_P}{}_{\ \ell_0}^{\ell_1} , 
     {S^\prime}_{i_1 \cdots i_P}{}_{\ \ell_0}^{\ell_1}
  \in {\mathfrak{T}}_P^0 \otimes
         {\mathfrak{T}}_1^1$, similarly write    
\begin{equation*}
S_{i_1 \cdots i_P}{}_{\ell_0}{}^{\ell_1} 
\sim {S^\prime}_{i_1 \cdots i_P}{}_{\ell_0}{}^{\ell_1} \   
\modd (i_1 \cdots i_{P} ) \Leftrightarrow
S_{[i_1 \cdots i_P]}{}_{\ell_0}{}^{\ell_1} 
={S^\prime}_{[i_1 \cdots i_P]}{}_{\ell_0}{}^{\ell_1}.    
\end{equation*}

\bigskip
Let $LM=\{\gamma:S^1\to M\}$ be the space of smooth loops on $M$.
$T_\gamma LM = \{X_\gamma\in C^{\infty}(S^1,TM): X_\gamma(\theta) \in T_{\gamma(\theta)}M\}$ is the set of 
smooth vector fields along $\gamma.$

We recall the definition of the Wodzicki-Chern-Simons class on $LM$ \cite{MRT4}.
For the metric $h_\rho$, set
\begin{align}\label{eq:K}
K_{\nu \lambda_1 \cdots \lambda_{2k+1} }{}_{\kappa_0}{}^{{\kappa_0}^\prime}     
&= \bar{R}_{\lambda_1 \kappa_1 \nu} {}^{{\kappa_0}^\prime} 
\bar{R}_{\lambda_2 \lambda_3 \kappa_2}{}^{\kappa_1}
\bar{R}_{\lambda_4 \lambda_5 \kappa_3}{}^{\kappa_2}
\cdots 
\bar{R}_{\lambda_{2k} \lambda_{2k+1} \kappa_0}{}^{\kappa_k} \in \Lambda^{2k+2}(M)\otimes \mathfrak{T}^1_1(M),\notag\\
K_{\nu \lambda_1 \cdots \lambda_{2k+1} }
&=\tr( K_{\nu \lambda_1 \cdots \lambda_{2k+1} }{}_{\kappa_0}{}^{{\kappa_0}^\prime} )
=K_{\nu \lambda_1 \cdots \lambda_{2k+1} }{}_{\kappa_0}{}^{{\kappa_0}},\\
K_{\nu [{\lambda_1} \cdots {\lambda_{2k+1}}]} &\in  \Lambda^1(M)\otimes \Lambda^{2k+2}(M).\notag
\end{align}
We define $CS^W\in \Omega^{2k+1}(LM)$ at the loop $\gamma \in LM$  by
\begin{equation}\label{eq:18}
CS^W(\gamma ) (X_{\gamma, 1}, \cdots , X_{\gamma, 2k+1}) 
= \int_0^{2\pi} 
 K_{\nu [{\lambda_1} \cdots {\lambda_{2k+1}}]}
(\gamma (\theta)) \dot{\gamma}^\nu (\theta) 
X_{\gamma , 1}^{\lambda_1}(\theta) \cdots  
X_{\gamma , 2k+1}^{\lambda_{2k+1}}(\theta)\   d\theta.
\end{equation}

Associated to a $C^\infty$ map $a: S^1 \times M \to M$
are: (i) 
$a^D: S^1 \to {\rm Maps}(M,M)$,  
$a^D(\theta)(m) := a(\theta,m)$; (ii)  $a^L:M\to LM$, 
$a^L(m)(\theta):= a(\theta,m). $  If $\xi$ is regular, the action (\ref{action})  has $a^L(m)$
equal to the loop  $\gamma_m$ which is the circle fiber containing $m$ and starting at $m$.  From now on, we use this action.

We want to compute $(a^L)^* CS^W  \in \Omega^{2k+1}(M).$
We have
\begin{align*}
 ((a^L)^* CS^W)_m &= (a^L)^* CS^W(\partial_{x^{i_1}},\cdots,\partial_{x^{i_{2k+1}}})_m dx^{i_1}\wedge\cdots\wedge 
 dx^{i_{2k+1}}.
\end{align*}
In our case, $a$ is rotation in the $S^1$ fiber, and
\begin{align*}
 \MoveEqLeft{(a^L)^* CS^W 
 (\partial_{x^{i_1}} , \cdots , \partial_{x^{i_{2k+1}}})_m} \\
 &= \int_{0}^{2\pi} 
 K_{\nu [\lambda_1 \cdots \lambda_{2k+1}]}  
   (a(\theta,m))
  \frac{d}{d\theta}\biggl|_{_{\theta}}a^\nu(\theta,m)
  (a^L_* (\partial_{x^{i_1}}|_m))^{\lambda_1} \cdots
  (a^L_* ( \partial_{x^{i_{2k+1}}}|_m))^{\lambda_{2k+1}} 
 \ d\theta \notag \\
 &= 
\int_0^{2\pi} K_{\nu [\lambda_1 \cdots \lambda_{2k+1}]} 
   (a(\theta,m))   
   \frac{\partial a^\nu}{\partial x^j}\biggl|_{(\theta,m)}\xi^j|_m
   \frac{\partial a^{\lambda_1}}{\partial x^{i_1}}\biggl|_{(\theta,m)}\cdots \frac{\partial a^{\lambda_{i_{2k+1}}}}{\partial x^{i_{2k+1}}}
   \biggl|_{(\theta,m)}
 \ d\theta,\notag
\end{align*}
where $m=(x^1, \cdots, x^{2k+1} )$.  In the last line, we  use that $a$ is an action, {\em i.e.,} $a(\theta + \psi, m) = a(\theta,a(\psi,m))$, to compute $(d/d\theta) a^\nu.$

For the action (\ref{action}),
 $a^D: S^1 \to \ism(M,h)$, so by the tensorial nature of the curvature tensor, we have
\begin{equation*}
 K_{\nu \lambda_1 \cdots \lambda_{2k+1}} 
   (a(\theta, m))   
   \frac{\partial a^\nu}{\partial x^j}|_{(\theta,m)}\xi^j|_m
   \frac{\partial a^{\lambda_1}}{\partial x^{i_1}}\biggl|_{(m,\theta)}\cdots 
    \frac{\partial 
     a^{\lambda_{2k+1}}}{\partial x^{i_{2k+1}}}
   \biggl|_{(m,\theta)}
   =  K_{j i_1 \cdots i_{2k+1}} \xi^j (m), 
\end{equation*}
for all $\theta.$
Thus, we have
\begin{align}\label{eq:a*CSW}
 ( (a^L)^* CS^W)  (m)  &=  
\left(\int_0^{2\pi} 
K_{j [i_1 \cdots i_{2k+1}]} (m) \xi^j (m)  
d\theta\right) 
dx^{i_1}  \wedge \cdots \wedge dx^{i_{2k+1}}\\
&= 
K_{j [i_1 \cdots i_{2k+1}]} (m) \xi^j (m) 
dx^{i_1}  \wedge \cdots \wedge dx^{i_{2k+1}}.
 \nonumber
\end{align}

\medskip
\par
As explained at the beginning of this section, our goal is to show that $( (a^L)^* CS^W)  (m) \neq 0$ for $h_\rho.$ 
This follows from the following Proposition, where it is more convenient to use indices $i_1,\cdots, i_{4k+1}$.
\begin{prop}\label{prop:1}
Let dim$(M) =4k+1$.  For $K_{j [i_1\cdots i_{4k+1}]}$ given by (\ref{eq:K}), we have
\begin{equation*}
K_{j [i_1 \cdots i_{4k+1}]}  \xi^j 
=
-(1+\rho^2)^2
(a_{4k} \rho^{4k} +  a_{4k-2}\rho^{4k-2} 
  \cdots + a_0), 
\end{equation*}
\label{formula: a_{4k}}
for
\begin{align*}
    a_{4k} &= (-1)^{k} 4^k(4k+2)
 \cdot  \eta_{[i_1} \phi_{{i_2}{i_3}} \cdots 
    \phi_{i_{4k} i_{4k+1}]}\\
    &= (-1)^{k} 4^k(4k+2)
    \left( \eta\wedge (d\eta)^k\right)_{i_1\cdots i_{4k+1}},   
\end{align*}
in the notation of (\ref{wedge-product}).
\end{prop}

The proof is in Appendix \ref{app:prop1}. From now on, we write $\eta\wedge (d\eta)^k$ just as $\eta\wedge d\eta^k$.  

\begin{cor} \label{cor:a4k} Let dim$(M) =4k+1$. For $m\in M$ and $\rho \gg 0$, 
$((a^L)^*CS^W)(m)\neq 0$ for the metric $h_\rho.$
\end{cor}
\noindent {\it Proof of  Corollary \ref{cor:a4k}.} 
By Lem.~\ref{lem:one}, $a_{4k}\neq 0.$  For $\rho \gg 0$,
$K_{j [1 \cdots 4k+1]}  \xi^j \neq 0$ by Prop~\ref{prop:1}.  It follows from (\ref{eq:a*CSW}) that
$((a^L)^*CSW)(m)\neq 0$. \hfill $\Box$

\medskip
Using the Corollary, we can prove the main theorem:
\begin{thm}\label{mainthm} Let 
 $(M, g, \phi , \xi , \eta )$ be $(4k+1)$-dimensional 
 closed regular Sasakian manifold, where $k \ge 1$.  
 Then, there is a metric of the form
 $h_\rho = g + \rho^2 \eta \otimes \eta$
 such that $|\pi_1(\ism(M, h_\rho)| = \infty$
 for $\rho\in \R, \rho \gg 0.$  For $\rho \neq \rho',$ $h_\rho$ and $h_{\rho'}$ are not isometric.
 \end{thm}

\noindent {\it Proof of Theorem \ref{mainthm}.} By 
 Cor.~\ref{cor:a4k}, $(a^L)^*CS^W$ for $h_\rho$ is pointwise nonzero for $\rho \gg 0.$ Since the rotation action (\ref{action}) is an action by isometries for $h_\rho$
 (Cor.~\ref{cor:hrho}), then
  \cite[Prop.~5.2]{MRT4} (as corrected in \cite[Prop.~2.2]{MRT5}) implies that $\pi_1(\ism(M, h_\rho))$ is infinite for $\rho \gg 0.$

  The curvature formula in Lem.~\ref{lem:3.3} shows that $h_\rho$ and $h_{\rho'}$ are not isometric for $\rho \gg 0.$
\hfill $\Box$

\section{Sasakian space forms} 
For the Sasakian space forms defined below,  the computations are more explicit.  
We start with the case of the sphere $S^{4k+1}$ with the standard metric 
$g_{std}$, and then treat the other Sasakian space forms. 

As in Ex.~\ref{exm:1}, $S^{4k+1}$ has the Hopf fibration 
\begin{equation}\label{hop} 
 S^1 \to S^{4k+1} \to {\mathbb{CP}}^{2k}.
 \end{equation}
Let $\xi$ be the unit vector field on $S^{4k+1}$
tangent to the $S^1$ fiber, 
and let $\eta$ be its dual $1$-form.   Then 
$(S^{4k+1}, g_{std}, \phi = \nabla \xi, \xi, \eta)$ is
a Sasakian structure on $S^{4k+1}$.

We can refine Prop.~\ref{prop:1}
for $S^{4k+1}$. 

\begin{prop}
\label{sphere-N}
Let $(S^{4k+1},g_{std}, \phi, \xi, \eta)$ 
be as above.
Then 
in the notation of Prop.~\ref{prop:1} 
\begin{equation} \label{eq:39}
K_{j[i_1 \cdots i_{4k+1}]} \xi^j 
=(1+\rho^2)^2  a_{4k} \rho^{4k}, 
\end{equation}
with
\begin{equation*}
a_{4k} =
(-1)^k 4^k(4k+2)(\eta \wedge d\eta^{2k})_{i_1 \cdots i_{4k+1}}.    
\end{equation*}
\end{prop}

Assuming the Proposition, we obtain:
\begin{thm}\label{thm:sphere}
Let $S^{4k+1}$ be $(4k+1)$-dimensional sphere
with the standard metric $g_{std}$.  
Then, we have
\begin{equation}
 |\pi_1 (\ism(S^{4k+1}, h_\rho))|= \infty   
\end{equation}
for the metric 
$h_\rho = g_{std} +\rho^2 \eta \otimes \eta$, 
$\rho \neq 0$.  
\end{thm}
Note the contrast to the standard metric at $\rho = 0$, where $\pi_1(\ism(S^{4k+1}, g_{std})) = \Z_2.$

The proof of the Theorem is almost the same as for Thm.~\ref{mainthm}: since 
$K_{j[i_1 \cdots i_{4k+1}]} \xi^j \neq 0$ in (\ref{eq:39}) and
the coefficients $a_r$
 of $\rho^{r}$ in  Prop.~\ref{formula: a_{4k}}  vanish, 
 $(a^L)^* CS^W \neq 0$ for any $\rho\neq 0.$  Then we follow the proof of Thm.~\ref{mainthm}.

\medskip
Prop.~\ref{sphere-N} will follow from the next Lemma. 
\begin{lem} \label{lem:5} For $(S^{4k+1}, h_\rho)$, in the notation of (\ref{123}) we have
\begin{align*}
(i) &\qquad 
(({\rm R}1)\otimes ({\rm R} {\rm E}))_{[i_2 i_3 i_4 i_5] \ell_3}{}^{\ell_1}   
=0;
\\
(ii) &\qquad 
(({\rm R}{\rm E})\otimes ({\rm R}1))_{[i_2 i_3 i_4 i_5] \ell_3}{}^{\ell_1}   
=0;
\end{align*}
for ${\rm E}=1,2,3.$
\end{lem}
\begin{proof} The proof is in Appendix \ref{app:S6}.
\end{proof}

\medskip
\noindent {\em Proof of Prop.~\ref{sphere-N}.} By Lem.~\ref{prop:RE-R3}, (\ref{eq:30}), and Lem.~\ref{lem:5},  the lowest power of $\rho$ in $K_{j[i_1 \cdots i_{4k+1}]} \xi^j$ is $\rho^{4k}.$ \hfill $\Box$
\bigskip

We can generalize this result  to Sasakian space forms. 
Recall that a Sasakian manifold $(M, \phi, \xi, \eta)$ is called 
a {\it Sasakian space form}
if all the $\phi$-sectional curvatures $K(X \wedge \phi X)$
are equal to a constant $c$, where $X$ is a unit vector field 
orthogonal to $\xi$, $\{X, \phi X\}$ span a two-plane section of $TM$, 
and $K(X \wedge \phi X)$ denotes the sectional curvature.
In this case, the
Riemanian curvature tensor of $M$ is given by
\begin{align}
\label{sasakian-space-form}
R_{kji}{}^h 
&= \frac{c+3}{4} (g_{ji} \delta_k{}^h -g_{ki}\delta_j{}^h )
\\
&\quad + \frac{c-1}{4} 
(\eta_i\eta_k \delta_j{}^h -\eta_j \eta_i \delta_k{}^h 
  +g_{ki}\eta_j \xi^h -g_{ji}\eta_k \xi^h
  +\phi_{ji} \phi_k{}^h 
    -\phi_{ki}\phi_j{}^h
   -2\phi_{kj} \phi_i{}^h ) 
   \notag 
\end{align}
\cite[Thm.~7.19]{Blair}.

As with ordinary space forms, Sasakian space forms fall into three classes \cite[Thm.~7.20]{Blair}, 
\cite{Tanno}
depending on whether
$c>-3$, $c=-3$, or $c<-3$.  
For  dim$( M) =4k+1$,
we have the following after scaling:

\begin{itemize}
    \item $c+3=4$.  $M$ is locally isometric
    to $(4k+1)$-dimensional sphere of constant sectional curvature $1$.
    \item $c+3=0$. $M$ is locally isometric to 
    the Heisenberg nilpotent Lie group $N = \R^{4k}\times \R$, with $g = \frac{1}{4}\sum_{i=1}^{2k}\left( (dx^i)^2 +  (dy^i)^2\right) + \eta\otimes\eta,$ with $\eta = \frac{1}{2}\left(dx^{4k+1} - \sum_{i=1}^{2k} y^idx^i\right).$
    \item$c+3=-1$. $M$ is locally isometric to 
    $\overline B = B\times \R$, where 
$(B, J, G)$ is a simply connected bounded domain in $\C^n$ with K\"ahler structure $J$ and metric $G$, the curvature of $G$ equals $\Omega = d\omega $, and  $\eta = \pi^*\omega + dt, g = \pi^*G + \eta\otimes \eta$ on $\overline B$ for the projection $\pi:\overline B \to B.$
\end{itemize}

Let $(M, g, \phi, \xi, \eta)$ be a $(4k+1)$-dimensional
Sasakian space form.  Its curvature tensor
(\ref{sasakian-space-form}) can be rewritten as
\begin{align}
 R_{kji}{}^h 
 &= \frac{c+3}{4} (g_{ji} \delta_k{}^h -g_{ki}\delta_j{}^h)
\notag  \\
&\qquad -\frac{c-1}{4}  (\phi_{ki} \phi_j{}^h -\phi_{ji}\phi_k{}^h
   +2\phi_{kj} \phi_i{}^h +2\eta_i \eta_k \delta_j{}^h
   -2\eta_j \eta_i \delta_k{}^h + g_{ki}\eta_j \xi^h
    -g_{ji}\eta_k \xi^h )
\notag \\
&\qquad + \frac{c-1}{4}(3\eta_i \eta_k \delta_j{}^h 
          -3\eta_j \eta_i \delta_k{}^h 
          +2g_{ki}\eta_j \xi^h -2g_{ij}\eta_k \xi^h) 
\notag \\
&= (\sr 1)_{kji}{}^h  + (\sr 2)_{kji}{}^h +(\sr 3)_{kji}{}^h,  
\notag 
\end{align}
where 
\begin{align*}
 (\sr 1)_{kji}{}^h & 
 = \frac{c+3}{4} (g_{ji} \delta_k{}^h -g_{ki}\delta_j{}^h),
 \\
 (\sr 2)_{kji}{}^h 
 &=-\frac{c-1}{4} (\phi_{ki} \phi_j{}^h 
 -\phi_{ji}\phi_k{}^h
   +2\phi_{kj} \phi_i{}^h +2\eta_i \eta_k \delta_j{}^h
   -2\eta_j \eta_i \delta_k{}^h + g_{ki}\eta_j \xi^h
    -g_{ji}\eta_k \xi^h ),
\notag \\
(\sr 3)_{kji}{}^h 
&= \frac{c-1}{4}
(3\eta_i \eta_k \delta_j{}^h 
          -3\eta_j \eta_i \delta_k{}^h 
          +2g_{ki}\eta_j \xi^h -2g_{ij}\eta_k \xi^h). 
          \notag 
\end{align*}

\begin{lem}\label{lem:6.2} 
\begin{equation*}
(\sr \alpha)_{i_2 i_3 \ell_2}{}^{\ell_1} 
       \times ({\rm RE})_{i_4 i_5 \ell_3}{}^{\ell_2} \sim 0 , 
\quad 
({\rm RE})_{i_2 i_3 \ell_2}{}^{\ell_1}  
\times
       (\sr\alpha)_{i_4 i_5 \ell_3 }{}^{\ell_2} \sim 0, 
\end{equation*}
for $\alpha = 1, {\rm E} = 1, 2, 3,$ and for ${\rm E} = 3, \alpha = 1, 2, 3.$ 
\end{lem}

\begin{proof} For $\alpha = 1,$ $(\sr 1) = (\rr 1),$ so the proof of Lem.~\ref{lem:5} carries over.
For ${\rm E} = 3,$  each term in  $(\sr 3)_{kji}{}^h$ contains either two $\eta$'s or an
$\eta\xi$ pair, so these terms are similar to $0$. \end{proof}

We note that 
\begin{align*}
\bar{R}_{kji}{}^h 
({\rm R}1)_{kji}{}^h
+(\sr 2^\prime)_{kji}{}^h
+(\sr 3)_{kji}{}^h 
 +({\rm R}3)_{kji}{}^h,   
\end{align*}
where
\begin{equation*}
(\sr 2^\prime)_{kji}{}^h    
=
-\left( \frac{c-1}{4}+ \rho^2\right)
(\phi_{ki} \phi_j{}^h 
  -\phi_{ji}\phi_k{}^h
   +2\phi_{kj} \phi_i{}^h 
     +2\eta_i \eta_k \delta_j^h
   -2\eta_j \eta_i \delta_k{}^h + g_{ki}\eta_j \xi^h
    -g_{ji}\eta_k \xi^h ).
\end{equation*}
Thus up to a scale factor $(\sr 2^\prime)_{kji}{}^h  = (\rr 2)_{kji}{}^h$ by (\ref{123}).  We can repeat the proof of Prop.~\ref{prop:1} by
replacing the power $\rho^{4k} = (\rho^2)^{2k}$ in Prop.~\ref{prop:1}
with $\frac{c-1}{4} +\rho^2$. We obtain: 
\begin{prop}
\label{sasakian-space-form-N}
Let $(M,g_0, \phi, \xi, \eta)$ be a 
Sasakian space form 
 as above.
Then 
in the notation of Prop.~\ref{prop:1} 
\begin{equation*} \label{eq:39a}
K_{j[i_1 \cdots i_{4k+1}]} \xi^j 
=(1+\rho^2)^2 a_{4k} \left(\frac{c-1}{4} +\rho^2\right)^{2k},  
\end{equation*}
with
\begin{equation*}
a_{4k} = (-1)^k 4k(4k+2)
(\eta \wedge d\eta^{2k})_{i_1 \cdots i_{4k+1}}.    
\end{equation*}
\end{prop}

As before, we conclude the following:

\begin{thm}
 Let $(M,g,\phi, \xi, \eta)$  be $(4k+1)$-dimensional
 closed orientable Sasakian space form. Then 
 $|\pi_1(\ism(M,g_\rho))|= \infty$ for 
 (i) $c+3>4, \rho \neq 0$, 
 (ii) $c+3=0, \rho \neq \pm 1,$ 
 or 
 (iii) $c+3< 0, \rho^2 \neq -\frac{c-1}{4}$.  
 \end{thm}

\begin{rem}
It follows from \cite[Note 7.5]{bk}  that $|\pi(\ism (M,g_0))|=\infty$   
for the Heisenberg infra-nilmanifold $M = N/\Gamma$
({\em i.e.,} a Sasakian space form with $c+3=0$).
This follows also from our result.  
\end{rem}

\section{Strict contact diffeomorhisms of $S^{4k+1}$}

We consider again the 
sphere $S^{4k+1}$ with the standard metric $g_{std}$,
and  Sasakian structure $(g_{st}, \phi, \eta, \xi)$.
The Reeb vector field $\xi$ of the Hopf fibering
\begin{equation*}
 \pi: S^1 \to S^{4k+1} \to \mathbb{CP}^{2k}    
\end{equation*}
is the unit tangent vector field along the $S^1$ fiber. It generates not only
an $S^1$-action via isometries of $g_{std}$, but also an $S^1$-action of strict contact diffeomorphisms of
$S^{4k+1}$ as a contact manifold.  Here, we
 recall that a diffeormorphism $\phi$ of $S^{4k+1}$
is  strictly contact if 
$\phi^* \eta = \eta.$
We denote the group of strict contact 
diffeomorphisms of $(S^{4k+1},\eta)$ by 
$\diff_{\eta, str} (S^{4k+1},\eta).$  We note that 
this group is infinite dimensional.

Then we have the following:
\begin{thm}
\label{thm:7.1}
\begin{equation*}
 |\pi_1(\diff_{\eta, str} (S^{4k+1}, \eta))| =
\infty.    
\end{equation*}
\end{thm}

The proof is similar to those in other sections, and will be given below.
As  in \S4, we take the family of 
metrics $h=g_\rho$ defined by
$h_{ij} = g_{ij} +\rho^2 \eta_i \eta_j $ for $g = g_{std}.$.
We can explicitly compute the integrand 
$K_{\nu[\lambda_1 \lambda_2 \cdots \lambda_{4k+1}]}$
in the definition of  $CS^W$ (\ref{eq:18}).

\begin{prop}
\label{prop:7.1}
\begin{equation*}
K_{\nu [\lambda_1 \lambda_2 \cdots \lambda_{4k+1}]}
=C \rho^{4k}  g_{\nu [\lambda_1}
 ((d\eta)^{2k})_{\lambda_1 \lambda_2 \cdots \lambda_{4k+1}]}, 
\end{equation*}
for some nonzero constant $C$.   
\end{prop}

The proof of Prop.~\ref{prop:7.1} is given after the next two Lemmas.

Recall that
\begin{equation*}
K_{\nu\lambda_1 \lambda_2 \cdots \lambda_{4k+1}}
= {\bar{R}_{\lambda_1 \kappa_1 \nu}}{}^{\kappa_0}
 {\bar{R}_{\lambda_2 \lambda_3 \kappa_2}}{}^{\kappa_1}
 \cdots 
{\bar{R}_{\lambda_{4k} \lambda_{4k+1} \kappa_{0}}}{}^{\kappa_{2k}}.
\end{equation*}
Using (\ref{eq:14}) and (\ref{constant-curvature form}), we  have n 
\begin{equation*}
\bar{R}_{kji}^{h}  
=(\crr 1)_{kji}{}^h 
-\rho^2 (\crr 2)_{kji}{}^h
-\rho^2(\crr 3)_{kji}{}^h 
-(\rho^2 + \rho^4) (\crr4)_{kji}{}^h, 
\end{equation*}
where 
\begin{align*}
(\crr1)_{kji}{}^h 
  &= g_{ji} \delta_k^h -g_{ki}\delta_j^h,   
\ (\crr2)_{kji}{}^h 
=\phi_{ki}\phi_j{}^h -\phi_{k}{}^h
  +2\phi_{ki}\phi_i{}^h,
  \notag \\
(\crr3)_{kji}{}^h 
&= g_{ki}\eta_j \xi^h -g_{ji}\eta_k \xi^h,
\ (\crr4)_{kji}{}^h 
= \eta_k\eta_j \delta_j{}^h 
   -\eta_j \eta_i \delta_k{}^h.
\notag 
\end{align*}
  
\begin{lem}
\label{CR}
 \begin{align*}
(\crr1)_{i_1 \kappa_1 j}{}^{\kappa_0}  \times
\bar{R}_{i_2 i_3 \kappa_2}{}^{\kappa_1} \cdots 
\bar{R}_{i_{4k} i_{4k+1} \kappa_0}{}^{\kappa_{2k}}
& 
\sim
-g_{\i_1 j}
\bar{R}_{i_2 i_3 \kappa_1}{}^{\kappa_0} \cdots 
\bar{R}_{i_{4k} i_{4k+1} \kappa_0}{}^{\kappa_{2k}}
\quad 
\modd (i_1 i_2 i_3 ).
\notag
 \end{align*}   
\end{lem}

\begin{proof}
\begin{align*}
\MoveEqLeft{
(\crr1)_{i_1 \kappa_1 j}{}^{\kappa_0}  \times
\bar{R}_{i_2 i_3 \kappa_2}{}^{\kappa_1} \cdots 
\bar{R}_{i_{4k} i_{4k+1} \kappa_0}{}^{\kappa_{2k}}   }\\
&= (g_{\kappa_1 j}\delta_{i_1}{}^{\kappa_0} 
  - g_{i_1 j}\delta_{\kappa_1}{}^{\kappa_0})
\bar{R}_{i_2 i_3 \kappa_2}{}^{\kappa_1} \cdots 
\bar{R}_{i_{4k} i_{4k+1} \kappa_0}{}^{\kappa_{2k}} 
\notag \\
&=\bar{R}_{i_2 i_3 \kappa_2}{}^{\kappa_1} \cdots 
\bar{R}_{i_{4k} i_{4k+1} i_1}{}^{\kappa_{2k}} 
-g_{i_1 j} 
\bar{R}_{i_2 i_3 \kappa_2}{}^{\kappa_0} \cdots 
\bar{R}_{i_{4k} i_{4k+1} \kappa_0}{}^{\kappa_{2k}}
\notag \\
&\sim 
-g_{i_1 j} 
\bar{R}_{i_2 i_3 \kappa_2}{}^{\kappa_0} \cdots 
\bar{R}_{i_{4k} i_{4k+1} \kappa_0}{}^{\kappa_{2k}},
\notag 
\end{align*}
because
$\bar{R}_{[i_{4k} i_{4k+1} i_1]}{}^{\kappa_{2k}}=0$.
\end{proof}

Set 
\begin{equation*}
T_{j i_1 i_2 \cdots i_{4k+1} }
= \bar{R}_{i_2 i_3 \kappa_2}{}^{\kappa_0} 
 \bar{R}_{i_4 i_5 \kappa_3}{}^{\kappa_2} 
\cdots 
\bar{R}_{i_{4k} i_{4k+1} \kappa_0}{}^{\kappa_{2k}}.
\end{equation*}

By Lemma \ref{CR}, we have
\begin{equation*}
K_{j [i_1 i_2 \cdots i_{4k+1}] }  
= -g_{j[i_1}T_{i_2 i_3 \cdots i_{4k+1}]}.
\end{equation*}
Many terms in $T_{j i_1 i_2 \cdots i_{4k+1} }$
vanish by direct computations as in Lemmas \ref{lem:5.4}, \ref{lem:6.2}.

\begin{lem}
$$(\crr{\rm E})_{i_2 i_3 j_1}{}^{j_0^\prime} \times 
(\crr{\rm F})_{i_4 i_5 j_2}{}^{j_1} \sim 0 \qquad \modd (i_1 i_2 i_3 i_4),
$$
for $E, F\in \{1,2,3,4\}$, $(E,F)\neq (2,2).$
\end{lem}

\noindent {\em Proof of Prop.~\ref{prop:7.1}.}  
We have
\begin{equation*}
T_{j i_1 i_2 \cdots i_{4k+1}}
\sim \rho^{4k}(\crr2)_{i_2 i_3 j_2}{}^{j_0}
(\crr2)_{i_2 i_3 j_2}{}^{j_0}
\cdots 
(\crr2)_{i_{4k} i_{4k+1} j_0}{}^{j_{2k}}\qquad \modd(i_1\cdots i_{4k+1}).
\end{equation*}
As in \S 5,  there is a nonzero constant $C$ such that  
\begin{equation*}
T_{j [i_1 i_2 \cdots i_{4k+1}]}
= C \rho^{4k}((d\eta)^{2k})_{i_2 i_3 \cdots i_{4k+1}}
\end{equation*}
${}$\hfill $\Box$
\medskip

\noindent {\em Proof of Thm.~\ref{thm:7.1}.}
By Proposition \ref{prop:7.1},  
$K_{j [i_1 i_2 \cdots i_{4k+1}] }\xi^j $ is given explicitly
in terms of $\eta$ and $d\eta$ (see Prop.~\ref{prop:1}).  Thus,
$K_{j [i_1 i_2 \cdots i_{4k+1}] }\xi^j $ is invariant under
 strict contact diffeomorphisms. For the $S^1$-action of the Reeb vector field $\xi$, the proof of Thm.~\ref{thm:5.1} can be directly adapted
 to $\diff_{\eta, str}(S^{4k+1}).$
 \hfill $\Box$

\begin{rem}
(i) By a slight modification of these computations, Thm.~\ref{thm:7.1} can be easily extended to the Sasakian space forms  of \S6.

(ii) Thm.~\ref{thm:7.1} has been obtained in \cite[Cor.~7]{casals}, with for $S^{4k+1}$ generalized to a closed contact manifold $M$ with closed Reeb orbits.  While both approaches prove the same element  $[a^D]\in \pi_1(\diff_{\eta, str}(M))$ has infinite order, \cite{casals} uses the map 
$a^D:S^1\to \diff_{\eta,str}(M)$ to pull back characteristic classes from $B\diff_{\eta,str}(M)$ to $BS^1.$
In contrast, we use the equivalent map $a^L:M\to LM, a^L(m)(\theta) = a^D(\theta)(m)$, to pull back a (nonclosed) Chern-Simons form from $LM$ to $M.$  It would be interesting to understand the relationship between these approaches better. In \cite{MRT6}, we reprove \cite[Cor.~7]{casals} using our approach. 
\end{rem}

\appendix

\section{Metric and curvature computations for the metric $h= h_\rho$}\label{app:curv}

\begin{customlem}{4.2}   For the metric $h = h_\rho,$
\begin{equation*}
  \bar{\Gamma}_{ij}^{k}  
    =  \Gamma_{ij}^{k}   
     - \rho^2 (\phi_i^k \eta_j + \phi_j^k \eta_i). 
\end{equation*}    
\end{customlem}

\begin{proof}
\begin{align*}
 \bar{\Gamma}_{ij}^{k} 
 &=
\frac{1}{2} h^{k\ell}
    (\partial_i g_{\ell j}+ \partial_j g_{i\ell} 
          -\partial_{\ell} g_{ij} ) 
+\frac{1}{2} \rho^2 h^{k\ell}   
( \partial_i (\eta_\ell \eta_j)
   +\partial_j (\eta_i \eta_\ell ) - \partial_\ell (\eta_i \eta_j))
 \\
&= 
\frac{1}{2} g^{k\ell} (\partial_i g_{\ell j} 
    + \partial_j g_{i\ell} - \partial_\ell g_{ij} ) 
+ \frac{1}{2} \alpha \xi^k \xi^\ell (\partial_i g_{\ell j} 
      +\partial_j g_{i\ell} -\partial_\ell g_{ij} ) 
\notag\\
&\qquad  +
\frac{1}{2} g^{k\ell} \rho^2 (\partial_i(\eta_\ell \eta_j) 
+ \partial_j(\eta_i \eta_\ell ) -\partial_\ell (\eta_i \eta_j))  
+
\frac{1}{2} \alpha \rho^2 \xi^k \xi^\ell  
   (\partial_i (\eta_\ell \eta_j) 
   + \partial_j (\eta_i \eta_\ell )
      - \partial_\ell (\eta_i \eta_j))  
      \notag\\
&:= \Gamma^k_{ij} +(b) + (c) + (d).
\notag 
\end{align*}

Note that
\begin{align*}
 \partial_i g_{\ell j} 
&=  \partial_i\langle\partial_{\ell},\partial _{j}\rangle
= \langle \nabla_{i} \partial_\ell, \partial_j\rangle + \langle \partial_\ell,\nabla_i\partial_j\rangle   
= 
  \Gamma_{i \ell}^{m} g_{mj} 
    + \Gamma_{ij}^{m} g_{\ell m}. 
    \notag
\end{align*}

Using 
$\Gamma_{ij}^k = \Gamma_{ji}^k$, we have 
\begin{align}
 (b) &= 
\frac{1}{2} \alpha \xi^k \xi^\ell (\partial_i g_{\ell j} 
      +\partial_j g_{i\ell} -\partial_\ell g_{ij} ) 
\notag \\
&= 
\frac{1}{2} \alpha \xi^k \xi^\ell
( 
\Gamma_{i \ell}^{m} g_{m j} 
+\Gamma_{ij}^{m} g_{\ell m} 
+\Gamma_{ji}^{m} g_{m\ell} 
+\Gamma_{j\ell}^{m} g_{i m} 
-\Gamma_{\ell i}^{m} g_{mj} 
-\Gamma_{\ell j}^{m} g_{i m} )  
\notag \\
&= \alpha \xi^k \xi^\ell 
\Gamma_{ij}^{m} g_{\ell m} 
= \alpha
\Gamma_{ij}^{m} \eta_m \xi^k, 
\notag
\end{align}
since $\xi^\ell g_{\ell m} = \eta_m$. 

Using (\ref{eq:con}) 
and $\partial_i \eta_j = \nabla_i \eta_j 
   + \Gamma_{ij}^k \eta_k $,
  we get 
\begin{align*}
(c) &= 
\frac{\rho^2}{2} g^{k\ell}
(\partial_i \eta_\ell\cdot \eta_j + \eta_\ell \cdot\partial_i\eta_j
+ \partial_j \eta_i \cdot\eta_\ell + \eta_i\cdot \partial_j \eta_\ell   
-\partial_\ell \eta_i \cdot\eta_j - \eta_i \cdot\partial_\ell \eta_j )
\notag \\
&=  \frac{\rho^2}{2}
g^{k\ell} 
( \nabla_i \eta_\ell\cdot \eta_j + \eta_\ell\cdot \nabla_i \eta_j
 + \nabla_j \eta_i \cdot\eta_\ell + \eta_i\cdot \nabla_j \eta_\ell   
-\nabla_\ell \eta_i\cdot \eta_j - \eta_i\cdot \nabla_\ell \eta_j 
\notag \\
&\qquad + \Gamma_{il}^m \eta_m \eta_j + \eta_\ell \Gamma_{ij}^m \eta_m
+ \Gamma_{ji}^m \eta_m \eta_\ell 
+ \eta_i \Gamma_{j\ell}^m \eta_m
-\Gamma_{\ell i}^m \eta_m \eta_j -\eta_i \Gamma_{\ell j}^m \eta_m )
\notag \\
&= \frac{\rho^2}{2} g^{k\ell} (\nabla_i\eta_\ell\cdot \eta_j 
-\nabla_\ell \eta_i \cdot\eta_j  + \eta_i\cdot \nabla_j \eta_\ell
-\eta_i\cdot \nabla_\ell \eta_j 
+2\Gamma_{ij}^{m} \eta_\ell\cdot \eta_m ) 
\notag \\
&= -\rho^2 (\phi_i{}^k \eta_j + \phi_j{}^k \eta_i ) 
+\rho^2 
\Gamma_{ij}^{m} \xi^k \eta_m, 
\notag 
\end{align*}
where (\ref{eq:refine}) is used in the last line.

Similarly,
\begin{align}
(d) &= 
\frac{1}{2} \alpha \rho^2 \xi^k \xi^\ell  
   (\partial_i (\eta_\ell \eta_j) 
   + \partial_j (\eta_i \eta_\ell )
      - \partial_\ell (\eta_i \eta_j))  
= 
-\frac{\alpha}{2} \xi^k \xi^\ell \rho^2 (2
 \phi_{i\ell} \eta_j + 2\phi_{j\ell} \eta_i 
    -2\Gamma_{ij}^{m} \eta_\ell \eta_m )
    \notag \\
&= \alpha \rho^2 
\Gamma_{ij}^{m} \xi^k \eta_m, 
\notag 
\end{align}
by (\ref{eq:one}).

Thus, 
\begin{equation*}
(b)+(c)+(d) = -\rho^2(\phi_i^k \eta_j + \phi_j^k \eta_i )
+ (\alpha +\rho^2 +\alpha\rho^2) 
    \Gamma_{ij}^{m} \xi^k \eta_m
=- \rho^2 (\phi_i{}^k \eta_j + \phi_j{}^k \eta_i ). 
\end{equation*}
This gives the lemma.
\end{proof}

\begin{customlem}{4.3}\label{lem:3.3a}
Let $\bar{R}_{kji}{}^h$ and $R_{kji}{}^h$
be the curvature tensor of the
Riemannian metrics $h_\rho$ and $g$ respectively.  
Then
\begin{align}\label{eq:14}
\bar{R}_{kji}{}^h &= R_{kji}{}^h 
 - \rho^2 (\phi_{ki} \phi_j{}^h -\phi_k{}^h \phi_{ji}
+  2 \phi_{kj} \phi_i{}^h 
+2\eta_k \eta_i \delta_j{}^h -2\eta_j \eta_i \delta_k{}^h 
+g_{ki} \eta_j \xi^h -g_{ji}\eta_k \xi^h ) \notag
\\
&\qquad - \rho^4(\eta_k \eta_i\delta_j{}^h  
- \eta_j \eta_i\delta_k{}^h ). 
\end{align}
\end{customlem}

\begin{proof}

\begin{align*}
\bar{R}_{kji}{}^{h} 
&= 
\partial_k \bar{\Gamma}_{ji}^{h}  
-\partial_j \bar{\Gamma}_{ki}^{h} 
+\bar\Gamma_{k \ell}^{h}
\bar\Gamma_{j i}^{\ell}
-\bar\Gamma_{j \ell}^{h}
\bar\Gamma_{ki}^{\ell}
\\
&=
\partial_k( \Gamma_{ji}^{h}
    -\rho^2 (\phi_j{}^h \eta_i + \phi_i{}^h \eta_j)  )
-\partial_j( \Gamma_{ki}^{h}  
    -\rho^2 (\phi_k{}^h \eta_i + \phi_i{}^h \eta_k)  ) 
\notag\\
&\qquad + (\Gamma_{k \ell}^{h}   
     -\rho^2 (\phi_k{}^h \eta_\ell+ \phi_\ell{}^h \eta_k ) )
    ( \Gamma_{ji}^{\ell}   
     -\rho^2 (\phi_j{}^\ell \eta_i
             + \phi_i{}^\ell \eta_j ) )
\notag \\
&\qquad - (\Gamma_{j \ell}^{h}   
     -\rho^2 (\phi_j{}^h \eta_\ell 
           +\phi_\ell{}^h \eta_j ) )
    ( \Gamma_{ki}^{\ell}   
     -\rho^2 (\phi_k{}^\ell \eta_i
               + \phi_i{}^\ell \eta_k ) )
 \notag 
\end{align*}
Therefore,
\begin{align*}
\bar{R}_{kji}{}^h &=
R_{kji}{}^h  
- \rho^2 (\partial_k(\phi_j{}^h \eta_i +\phi_i{}^h \eta_j )
  -\partial_j(\phi_k{}^h \eta_i +\phi_i{}^h \eta_k ))
\notag \\
&\qquad + 
\rho^2 ( \Gamma_{k\ell}^{h} 
(\phi_j^\ell\eta_i + \phi_i{}^\ell \eta_j) 
+ \Gamma_{ji}^{\ell}
(\phi_k{}^h\eta_\ell + \phi_\ell{}^h \eta_k) 
)
\notag \\
&\qquad -\rho^2 ( 
\Gamma_{j\ell}^{h} 
    (\phi_k^\ell\eta_i + \phi_i^\ell \eta_k) 
+\Gamma_{ki}^{\ell}  
     (\phi_j^h \eta_\ell + \phi_\ell^h \eta_j) 
) 
\notag \\
&\qquad + \rho^4 
(\phi_k{}^h \eta_\ell + \phi_\ell{}^h \eta_k) 
(\phi_j{}^\ell \eta_i + \phi_i{}^\ell \eta_j) 
-\rho^4
(\phi_j{}^h \eta_\ell + \phi_\ell{}^h \eta_j) 
(\phi_k{}^\ell \eta_i + \phi_i{}^\ell \eta_k) 
\notag \\
&= R_{kji}{}^h + (\rho_2) + (\rho_4) 
\notag 
\end{align*}  
where $(\rho_a)$ is the $\rho^{a}$-term of 
$\bar{R}_{kji}{}^h$. We have 
\begin{align}
(\rho_2) 
&= -[\rho^2 (\partial_k(\phi_j{}^h \eta_i 
       +\phi_i{}^h \eta_j )
  -\partial_j(\phi_k{}^h \eta_i +\phi_i{}^h \eta_k ))
\notag \\
&\qquad  +
\rho^2 ( 
\Gamma_{k\ell}^{h} 
(\phi_j{}^\ell\eta_i + \phi_i{}^\ell \eta_j) 
+ \Gamma_{ji}^{\ell}
(\phi_k{}^h\eta_\ell + 
           \phi_\ell{}^h \eta_k) 
)
\notag \\
&\qquad  -\rho^2 ( 
\Gamma_{j \ell}^{h}
(\phi_k{}^\ell\eta_i + \phi_i{}^\ell \eta_k) 
+ \Gamma_{ki}^{\ell}
     (\phi_j{}^h \eta_\ell + \phi_\ell{}^h \eta_j) 
) ]
\notag \\
&= 
-\rho^2 [
(\nabla_k \phi_j{}^h \cdot\eta_i + \phi_j{}^h\cdot \nabla_k\eta_i
+ \nabla_k \phi_i{}^h \cdot\eta_j 
           + \phi_i{}^h \cdot\nabla_k \eta_j ) 
\notag \\
& \qquad  
   -\rho^2(\nabla_j \phi_k{}^h \cdot\eta_i 
        + \phi_k{}^h \cdot\nabla_j\eta_i
+ \nabla_j \phi_i{}^h\cdot \eta_k 
       + \phi_i{}^h\cdot \nabla_j \eta_k )],
\notag 
\end{align}
 where we use
$$\nabla_k\eta_i := (\nabla_k\eta)_i = \partial_k\eta_i -\Gamma^\ell_{ki}\eta_\ell,
\ \nabla_k\phi_j{}^h := (\nabla_k\phi)_j{}^h = \partial_k \phi_j{}^h + \phi_j{}^\ell\Gamma^h_{k\ell} - \phi_\ell{}^h \Gamma^\ell_{kj}.$$    
Using (\ref{(A.2)}), (\ref{eq:refine}),  we have 
\begin{equation*}
 (\rho_2) = -\rho^2\left(
 g_{ki}\eta_j \xi^h -g_{ji}\eta_k \xi^h
 -2\delta_k{}^h \eta_j \eta_i + 2 \delta_j{}^h \eta_k \eta_i
 +\phi_j{}^h \phi_{ki} +2\phi_i{}^h \phi_{kj}
 -\phi_k{}^h \phi_{ji} \right).
\end{equation*}
For the $\rho_4$-term, 
using (\ref{eq:one}), we easily get
\begin{align*}
 (\rho_4)& =  \rho^4 
(\phi_k{}^h \eta_\ell + \phi_\ell{}^h \eta_k) 
(\phi_j{}^\ell \eta_i + \phi_i{}^\ell \eta_j)
-\rho^4
(\phi_j{}^h \eta_\ell + \phi_\ell{}^h \eta_j) 
(\phi_k{}^\ell \eta_i + \phi_i{}^\ell \eta_k) 
\notag \\
&=
\rho^4(-\delta_j{}^h \eta_k \eta_i + \delta_k{}^h \eta_j \eta_i).
\end{align*}
Thus, we have the lemma. 
\end{proof}

\section{Proof of Proposition \ref{prop:1}}\label{app:prop1}

This Appendix gives the proof of Prop.~\ref{prop:1}.  As an outline, by (\ref{eq:quick2}),
$K_{j [i_1 \cdots i_{4k+1}]}  \xi^j = -(1+\rho^2) \eta_{[i_1} 
  (\tr(\rk 1))_{i_2 \cdots i_{2k+1}]}$, where $(\rk 1)$ is defined in (\ref{eq:quick}).  By the notation defined in (\ref{eq:quick3}), (\ref{eq:quick4}), it suffices to prove (\ref{prop:N}).  This is proven in Corollary \ref{cor:quick5}(ii). 
\medskip

\begin{customprop}{5.1}
Let dim$(M) =4k+1$.  For $K_{j [i_1\cdots i_{4k+1}]}$ given by (\ref{eq:K}), we have
\begin{equation*}
K_{j [i_1 \cdots i_{4k+1}]}  \xi^j 
=
-(1+\rho^2)^2
(a_{4k} \rho^{4k} +  a_{4k-2}\rho^{4k-2} 
  \cdots + a_0), 
\end{equation*}
%\label{formula: a_{4k}}
for
\begin{align*}
    a_{4k} &= (-1)^{k} 4^k(4k+2)
 \cdot  \eta_{[i_1} \phi_{{i_2}{i_3}} \cdots 
    \phi_{i_{4k} i_{4k+1}]}\\
    &= (-1)^{k} 4^k(4k+2)
    \left( \eta\wedge (d\eta)^k\right)_{i_1\cdots i_{4k+1}},   
\end{align*}
in the notation of (\ref{wedge-product}).
\end{customprop}

\noindent {\em Proof.}
We first work in an arbitrary odd dimension 
$2k+1.$
Recall that  
\begin{equation}\label{eq:111}
   K_{j i_1 \cdots i_{2k+1}}
 {}_{\ell^\prime_0}{}^{\ell_0}\xi^j \\
= 
\bar{R}_{i_1 \ell_1 j}{}^{\ell_0} 
\bar{R}_{i_2 i_3 \ell_2}{}^{\ell_1} 
\bar{R}_{i_4 i_5 \ell_3}{}^{\ell_2} 
\cdots
\bar{R}_{i_{2k} i_{2k+1} \ell^\prime_0}{}^{\ell_k}\xi^j. 
\end{equation}
(\ref{(A.2)}) implies $\nabla_j \xi^k =  \-\phi_j^k$ \cite[p.~87]{Blair},
so 
\begin{align}\label{eq:a}
R_{kji}^{\ \ \ h}\xi^i
&= \nabla_k \nabla_j \xi^h 
-\nabla_j \nabla_k \xi^h
=-(\nabla_k \phi_j{}^h - \nabla_j\phi_k{}^h )
\notag \\
&=-[(g_{kj}\xi^h -\eta_j \delta_k{}^h)
-(g_{jk}\xi^h -\eta_k \delta_j{}^h)]
 \\
&=-(\eta_k \delta_j{}^h -\eta_j \delta_k{}^h).\notag
\end{align}
By Lem.~\ref{lem:3.3}, we get
\begin{align*} 
\bar{R}_{kji}{}^h\xi^i&= R_{kji}{}^h \xi^i
 - \rho^2 (\phi_{ki} \phi_j{}^h -\phi_k{}^h \phi_{ji}
+  2 \phi_{kj} \phi_i{}^h 
+2\eta_k \eta_i \delta_j{}^h 
   -2\eta_j \eta_i \delta_k{}^h 
+g_{ki} \eta_j \xi^h -g_{ji}\eta_k \xi^h ) \xi^i \notag
\\
&\qquad - \rho^4(\eta_k \eta_i\delta_j^h  
- \eta_j \eta_i\delta_k{}^h )\xi^i 
    \\
   &=
-(\eta_k \delta_j{}^h -\eta_j \delta_k{}^h)
-2\rho^2 (\eta_k \delta_j{}^h -\eta_j \delta_k{}^h)
-\rho^4(\eta_k \delta_j{}^h - \eta_j \delta_k{}^h) \notag
\\
&=-(1+\rho^2)^2
(\eta_k\delta_j{}^h - \eta_j\delta_k{}^h), 
\notag
\end{align*}
where we use (\ref{eq:a}) for the $\rho^0$ term, (\ref{eq:one}) and (\ref{eq:two}) for the $\rho^2$ term, and (\ref{eq:one}) for the $\rho^4$ term.
Then 
\begin{align*} 
K_{j 
 i_1 \cdots i_{2k+1}} 
    {}_{\ell^\prime_0}{}^{\ell_0} 
 &=\xi^j \bar{R}_{i_1 \ell_1 j}{}^{\ell_0} 
\bar{R}_{i_2 i_3 \ell_2}{}^{\ell_1}
\bar{R}_{i_4 i_5 \ell_3}{}^{\ell_2}
\cdots
\bar{R}_{i_{2k} i_{2k+1} \ell^\prime_0}{}^{\ell_k}  \xi^j
\notag\\
&= 
-(1+\rho^2)^2
(\eta_{i_1} \delta_{\ell_1}{}^{\ell_0} 
   -\eta_{\ell_1} \delta_{i_1}{}^{\ell_0}) 
\bar{R}_{i_2 i_3 \ell_2}{}^{\ell_1}
\bar{R}_{i_4 i_5 \ell_3}{}^{\ell_2}
\cdots
\bar{R}_{i_{2k} i_{2k+1} \ell^\prime_0}{}^{\ell_k} 
 \\
&:= (\rk 1)_{i_1 \cdots i_{2k+1}}
   {}_{\ell^\prime_0}{}^{\ell_0} 
+(\rk 2)_{i_1 \cdots i_{2k+1}}
     {}_{\ell^\prime_0}{}^{\ell_0},
\notag 
\end{align*}
where  
\begin{align}\label{eq:quick}
(\rk 1)_{i_1 \cdots i_{2k+1}} 
  {}_{\ell^\prime_0}{}^{\ell_0}
&= 
-(1+\rho^2)^2
(\eta_{i_1} \delta_{\ell_1}^{\ \ell_0} )
\bar{R}_{i_2 i_3 \ell_2}{}^{\ell_1}
\bar{R}_{i_4 i_5 \ell_3}{}^{\ell_2}
\cdots
\bar{R}_{i_{2k} i_{2k+1} \ell^\prime_0}{}^{\ell_k},  
\\
(\rk 2)_{i_1 \cdots i_{2k+1}} 
{}_{\ell^\prime_0}{}^{\ell_0}
&= 
(1+\rho^2)^2
(\eta_{\ell_1} \delta_{i_1}^{\ \ell_0}) 
\bar{R}_{i_2 i_3 \ell_2}{}^{\ell_1}
\bar{R}_{i_4 i_5 \ell_3}{}^{\ell_2}
\cdots
\bar{R}_{i_{2k} i_{2k+1} \ell^\prime_0}{}^{\ell_k}.
\notag 
\end{align}
As in (\ref{trace}), we have 
\begin{align*}
\tr(\rk 1)_{i_1 \cdots i_{2k+1}} 
&= 
 -(1+\rho^2)^2
\eta_{i_1}\bar{R}_{i_2 i_3 \ell_2}{}^{\ell_1}
\bar{R}_{i_4 i_5 \ell_3}{}^{\ell_2}
\cdots
\bar{R}_{i_{2k} i_{2k+1} \ell_1}{}^{\ell_k},
\\
\tr(\rk 2)_{i_1 \cdots i_{2k+1}} 
&= 
(1+\rho^2)^2
(\eta_{\ell_1}) 
\bar{R}_{i_2 i_3 \ell_2}{}^{\ell_1}
\bar{R}_{i_4 i_5 \ell_3}{}^{\ell_2}
\cdots
\bar{R}_{i_{2k} i_{2k+1} i_1}{}^{\ell_k}.
\end{align*}

We continue the proof of Prop.~\ref{prop:1} with a series of lemmas. 
Since the Proposition is concerned with the skew-symmetrization of $K\cdot\xi$, 
we consider the skew-symmetric parts of
$\tr(\rk 1)_{i_1 \cdots i_{2k+1}}$ and $\tr(\rk 2)_{i_1 \cdots i_{2k+1}}$. The following proof is the first instance of our repeated use of the Bianchi identity.

\begin{lem}
\begin{equation}
\tr(\rk 2)_{[i_1 \cdots i_{2k+1}]} =0.        
\end{equation}{}
\end{lem}

\begin{proof}
Set 
\begin{equation*}
(\rk 2\hy 1)_{i_2 \cdots i_{2k-1}}{}_{\ell_{k}}  
= \eta_{\ell_1}
\bar{R}_{i_2 i_3 \ell_2}{}^{\ell_1}
\bar{R}_{i_4 i_5 \ell_3}{}^{\ell_2}
\cdots
\bar{R}_{i_{2k-2} i_{2k-1} \ell_k}{}_{\ell_{k-1}},
\end{equation*}
and 
\begin{equation*}
\hat{\bar{R}}_{i_{2k-2} i_{2k-1} i_1}{}^{\ell_k} 
=\bar{R}_{[i_{2k-2} i_{2k-1} i_1]}{}^{\ell_k}.    
\end{equation*}
Then 
\begin{equation*}
\tr(\rk 2)_{[i_1 \cdots i_{2k+1}]}=   
\tr((\rk 2\hy 1) \wedge \bar{R})_{[i_1 \cdots i_{2k-1}]}. 
\end{equation*}
Note that 
${\bar{R}}_{[i_{2k-2} i_{2k-1} i_1]}{}^{\ell_k}=0$
by the Bianchi identity.  
Thus,  
$\tr(\rk 2)_{[i_1 \cdots i_{2k+1}]} =0. $
\end{proof}

Therefore,
\begin{equation}\label{eq:quick2}
 K_{j [i_1 \cdots i_{2k+1}]}\xi^j  
=-(1+\rho^2) \eta_{[i_1} 
  (\tr(\rk 1))_{i_2 \cdots i_{2k+1}]}.
\end{equation}

\bigskip
From now on, we focus on  $dim(M)=4k+1$. 
Let 
\begin{equation}   
\label{eq:quick3}
N_{i_2 i_3 \cdots i_{4k} i_{4k+1}}{}_{\ell_1}{}^{\ell_0} 
= \bar{R}_{i_2 i_3 \ell_2}{}^{\ell_0}
\bar{R}_{i_4 i_5 \ell_3}{}^{\ell_2}
\cdots 
\bar{R}_{i_{4k} i_{4k+1} \ell_1}{}^{\ell_{2k}}.
\end{equation}
Then 
\begin{align}\label{eq:quick4}
 (\rk 1)_{j[i_1 i_2 \cdots i_{4k+1}]} \xi^j 
& = -(1+\rho^2)^2 (\rk 1)_{[i_1 i_2 \cdots i_{4k+1}]}{}_{\ell}{}^{\ell}
= -(1+\rho^2)^2 \eta_{[i_1} 
 N_{i_2 \cdots \cdots i_{4k+1}]}{}_{\ell}{}^{\ell}
 \notag \\
& =-(1+\rho^2)^2 \sum_{a=1}^{4k+1} 
(-1)^a  \eta_{i_a} 
   N_{[i_1 \cdots \hat{i}_a \cdots i_{4k+1}]}{}_{\ell}{}^{\ell},
   \end{align}
where $\hat{i}_a$ means $i_a$ is omitted.
To get Proposition \ref{prop:1}, it suffices to prove 
\begin{align}\label{prop:N}
\tr(N)_{[i_2 \cdots i_{4k+1}]} &=  N_{[i_2 \cdots i_{4k+1}]}{}_{\ell}{}^{\ell}  
  = b_{4k}\rho^{4k} + \cdots + b_0,\\ 
  {\rm for}\ 
   b_{4k} &= (-1)^k 4^{k}(4k+2) (d\eta \wedge d\eta^k )_{i_2 \cdots i_{4k+1} }.  \notag
\end{align}

\bigskip
\par
We regroup 
$N_{i_2 i_3 \cdots i_{4k} i_{4k+1}}
      {{}_{\ell_0}{}^{\ell_1}}$
 as 
\begin{align} \label{eq:regroup}
&{N_{i_2 i_3 \cdots i_{4k} i_{4k+1}}}
   {}_{\ell_0}{}^{\ell_1} 
\\
&= A_{i_2 i_3 i_4 i_5 \ell_3}{}^{\ell_1}
\cdot
A_{i_6 i_7 i_8 i_9 \ell_5}{}^{\ell_3}
\cdots 
A_{i_{4k-6} i_{4k-5} i_{4k-4} i_{4k-3} \ell_{2k-1}}{}^{\ell_{2k-3}}
\cdot 
A_{i_{4k-2} i_{4k-1} i_{4k} i_{4k+1} \ell_{0}}
  {}^{\ell_{2k-1}},
  \notag
\end{align}
where
\begin{equation}\label{eq:A}
A_{i_{4a-2} i_{4a-1} i_{4a} i_{4a+1} \ell_{2a+1}}{}^{\ell_{2a-1}  }
:= 
\bar R_{i_{4a-2} i_{4a-1} \ell_{2a}}{}^{\ell_{2a-1}}
\cdot
\bar R_{i_{4a} i_{4a+1} \ell_{2a+1}}{}^{\ell_{2a}},\ a = 1,\ldots, k,
\end{equation} 
(with $\ell_{2k+1}$ replaced by $\ell_1$ on the left hand side of (\ref{eq:A}) if $a=k.$)
\bigskip

 {\it We emphasize that this regrouping is not possible if dim$(M) = 4k+3.$}

\bigskip

By Lem.~\ref{lem:3.3}, we have 
\begin{align*}
\bar{R}_{kji}{}^h  
&= ({\rm R}1)_{kji}{}^h + ({\rm R}2)_{kji}{}^h 
+({\rm R}3)_{kji}{}^h, 
\end{align*}
where
\begin{align}\label{123}
({\rm R}1)_{kji}{}^h &= R_{kji}{}^h, \notag \\
({\rm R}2)_{kji}{}^h &= 
-\rho^2 (\phi_{ki} \phi_j{}^h -\phi_k{}^h \phi_{ji}
+ 2 \phi_{kj} \phi_i{}^h 
+2\eta_k \eta_i \delta_j{}^h -2\eta_j \eta_i \delta_k{}^h + g_{ki}\eta_j\xi^h - g_{ji}\eta_k\xi^h), \\
({\rm R}3)_{kji}{}^h &= 
 \rho^4(\delta_k{}^h \eta_j \eta_i -\delta_j{}^h \eta_i \eta_k). \notag
\end{align}
In the notation of (\ref{123}) (and relaxing the sub/superscripts on $A$), we have
\begin{align} \label{eq:16}
A_{i_a i_{a+1} i_{a+2} i_{a+3} \ell_\beta}{}^{\ell_\alpha}  
& =
\left[({\rm R}1)_{i_{a} i_{a+1} \ell_\gamma}{}^{\ell_\alpha}
+
({\rm R}2)_{i_{a} i_{a+1} \ell_\gamma}{}^{\ell_\alpha}
+({\rm R}3)_{i_{a} i_{a+1} \ell_\gamma}{}^{\ell_\alpha}\right]\\
&\qquad \cdot 
\left[({\rm R}1)_{i_{a+2} i_{a+3} \ell_\beta}{}^{\ell_\gamma} +
({\rm R}2)_{i_{a+2} i_{a+3} \ell_\beta}{}^{\ell_\gamma}
+({\rm R}3)_{i_{a+2} i_{a+3} \ell_\beta}{}^{\ell_\gamma}\right].\nonumber
\end{align}

Thus, we will compute 
$({\rm RE})_{i_{a} i_{a+1} \ell_\gamma}{}^{\ell_\alpha}
\otimes 
({\rm RF})_{i_{a+2} i_{a+3} \ell_\beta}{}^{\ell_\gamma}$
for ${\rm E, F} =1,2,3$. For ease of notation,
we will just compute  $({\rm RE})_{i_{2} i_{3} \ell_2}{}^{\ell_1}
\otimes 
({\rm RF})_{i_{4} i_{5} \ell_3}{}^{\ell_2}$, 
as this easily extends to  the  general
case. 

\bigskip

The Bianchi identity extends to each component of the curvature tensor.
\begin{lem}
 \label{lem:RE-Bianchi}  
 \begin{equation}
 ({{\rm RE}})_{[kji]}{}^h =0, 
 \end{equation}
 for ${\rm E}=1,2,3$. 
\end{lem}
\begin{proof}
 The curvature tensor $\bar{R}_{kji}{}^h$ of $h_\rho$ has the Bianchi identity for    each value of $\rho.$  Letting $\rho$ vary in (\ref{eq:14}), we must have the Bianchi identity for each
 RE.
\end{proof}

We now  eliminate curvature terms involving $({\rm R}3)$.
\begin{lem} \label{prop:RE-R3}
\begin{equation}
\label{3 times 3}
 (({\rm RE}) \otimes({\rm R}3))_{[i_2i_{3} i_4 i_5]\ell_3}
 {}^{\ell_1}  
=
(({\rm R}3) \otimes({\rm RE}))_{[i_2i_{3} i_4 i_5]\ell_3}
 {}^{\ell_1} 
 =0,        
\end{equation}
for $E=1,2,3$.
\end{lem}
\begin{proof}  Ignoring powers of $\rho$, we have 
\begin{align*}
 ({\rm RE})_{i_2 i_3 \ell_2}{}^{\ell_1}  
 \cdot ({\rm R}3)_{i_4 i_5 \ell_3}{}^{\ell_2}
&=
({\rm RE})_{i_2 i_3\ell_2}{}^{\ell_1} 
(\delta_{i_4}{}^{\ell_2}
    \eta_{i_5}\eta_{\ell_2}
-\delta_{i_5}{}^{\ell_2}
    \eta_{i_4}\eta_{\ell_2})
\\
& 
=({\rm RE})_{i_2 i_3 i_4}{}^{\ell_1} \eta_{i_5}
-({\rm RE})_{i_2 i_3 i_5}{}^{\ell_1} \eta_{i_4}.
\notag 
\end{align*}
Using Lemma {\ref{lem:RE-Bianchi}}, we have
\begin{equation*}
(({\rm RE}) \otimes({\rm R}3))_{[i_2i_{3} i_4 i_5]}
{}_{\ell_3}{}^{\ell_1} 
 =  \eta_{[i_5}({\rm RE})_{i_2 i_3 i_4]}{}^{\ell_1} 
- \eta_{[i_4}({\rm R}E)_{i_2 i_3 i_5]}{}^{\ell_1} 
=0.
\end{equation*}
The second equation of \ref{3 times 3}
is  obtained similarly. 
\end{proof}

We recall the regrouping of 
$N_{i_2 i_3 \cdots i_{4k} i_{4k+1}}{}^{\ell_1}{}_{\ell_0}{}$ in (\ref{eq:regroup}) and that
each 
$A_{i_{4a-2} i_{4a-1} i_{4a} i_{4a+1} \ell_{2a+1}}
   {}^{\ell_{2a-1}  }$
is composed by $({\rm R}E)$'s.  

For $a > 4k$,
the terms in $\rho^{2a}$ 
of $N_{i_2 i_3 \cdots i_{4k} i_{4k+1}}{}_{\ell_0}{}^{\ell_1}$  should contain at least one $({\rm R}3)$.  

Thus, by Lemma \ref{prop:RE-R3}, we have 
\begin{align}\label{eq:30}
N_{[i_2 i_3 \cdots i_{4k} i_{4k+1}] }
 {}_{\ell_0} {}^{\ell_1}
&= 
\underbrace{(
({\rm R}2) \otimes \cdots \otimes ({\rm R}2)
)}_{\text{2k}}
{}_{[i_2 i_3 \cdots i_{4k} i_{4k+1}]} 
  {}_{\ell_0}{}^{\ell_1}  
  + l.o. \\
  &= \underbrace{(
({\rm R}2) \wedge \cdots \wedge ({\rm R}2)
)}_{\text{2k}}
{}_{i_2 i_3 \cdots i_{4k} i_{4k+1}} 
  {}_{\ell_0}{}^{\ell_1} 
  + l.o. \notag
\end{align}
where $l.o.$ denotes the lower order terms in $\rho^{2a}$, $2a \leq 4k-2$.
\medskip

We have
\begin{equation*}
 ({\rm R} 2)_{i_2 i_3 \ell_2}{}^{\ell_1}
 =
 -\rho^2 (({\rm R}2\hy 1)_{i_2 i_3 \ell_2}{}^{\ell_1}
  + ({\rm R}2\hy2)_{i_2 i_3 \ell_2}{}^{\ell_1} ), 
\end{equation*}
where
\begin{align*}
(\rr 2\hy 1)_{i_2 i_3 \ell_2}{}^{\ell_1}
&= \phi_{i_2 \ell_2}\phi_{i_3}{}^{\ell_1}
-\phi_{i_2}{}^{\ell_1}\phi_{i_3\ell_2}
+2 \phi_{i_2 i_3} \phi_{\ell_2}{}^{\ell_1},
\\
(\rr 2\hy2)_{i_2 i_3 \ell_2}{}^{\ell_1} 
&=2(\eta_{i_2} \delta_{i_3}{}^{\ell_1} 
      -\eta_{i_3}\delta_{i_2}{}^{\ell_1})
          \eta_{\ell_2}
 + (g_{i_2 \ell_2} \eta_{i_3} -g_{i_3 \ell_2}\eta_{i_2})     \xi^{\ell_1}.
     \notag 
\end{align*}

\begin{lem}\label{lem:5.4}
\begin{align}
\label{eq:R2-relation}
&((\rr 2\hy 1) \otimes (\rr 2\hy 2))_{[i_2 i_3 i_4 i_5] \ell_3}{}^{\ell_1}  =0,  
\notag\\
&((\rr 2\hy 2) \otimes (\rr 2\hy 1))_{[i_2 i_3 i_4 i_5] \ell_3}{}^{\ell_1}  =0,  
 \\
&((\rr 2\hy 2) \otimes (\rr 2\hy 2))_{[i_2 i_3 i_4 i_5] \ell_3}{}^{\ell_1}  =0. \notag 
\end{align}    
\end{lem}

\begin{proof}
A direct check gives
\begin{equation} \label{eq:(R2)bianki}
(\rr 2\hy 1)_{[kji]}{}^h =0 , 
(\rr 2\hy 2)_{[kji]}{}^h =0.    
\end{equation}

For the first equation of (\ref{eq:R2-relation}),
we get 
\begin{align*}
(\rr 2\hy 1)_{i_2 i_3  \ell_2}{}^{\ell_1} \times
(\rr 2\hy 2)_{i_4 i_5 \ell_3}{}^{\ell_2}
&=(\rr 2\hy 1)_{i_2 i_3 \ell_2}{}^{\ell_1}
\times 2(\eta_{i_4} \delta_{i_5}{}^{\ell_2}
  - \eta_{i_5} \delta_{i_4}{}^{\ell_2}) \eta_{\ell_3}
\notag \\
&=2(\eta_{i_4}(\rr 2\hy 1)_{i_2 i_3 i_5}{}^{\ell_1}
    - \eta_{i_5}(\rr 2\hy 1)_{i_2 i_3 i_4}{}^{\ell_1})\eta_{\ell_3},
    \notag 
\end{align*}
so
\begin{align*}
((\rr 2\hy 1)\otimes (\rr 2\hy 2)_{[i_2 i_3 i_4 i_5]}{}^{\ell_1} 
&= 2(\eta_{[i_4}(\rr 2\hy 1)_{i_2 i_3 i_5]}{}^{\ell_1}
    - \eta_{[i_5}(\rr 2\hy 1)_{i_2 i_3 i_4]}{}^{\ell_1})
=0. 
\notag 
\end{align*}

For the second equation of (\ref{eq:R2-relation}),
we have 
\begin{equation*}
(\rr 2\hy 2)_{i_2 i_3  \ell_2}{}^{\ell_1} \times
(\rr 2\hy 1)_{i_4 i_5 \ell_3}{}^{\ell_2} 
=
((\rr 2\hy 1)_{i_4 i_5 \ell_3 i_2 }\eta_{i_3} 
-(\rr 2\hy 1)_{i_4 i_5 \ell_3 i_3}\eta_{i_2}) \xi^{\ell_1}.
\end{equation*}
so
\begin{equation*}
((\rr 2\hy 2) \otimes (\rr 2\hy 1))_{[i_2 i_3 i_4 i_5] \ell_3}{}^{\ell_1}
= -(\rr 2\hy 1)_{[i_4 i_5 i_2]\ell_3}{}^{\ell_1}\eta_{[i_3]}
 + (\rr 2\hy 1)_{[i_4 i_5 i_3]\ell_2}{}^{\ell_1}\eta_{[i_2]}, 
\end{equation*}
with the alternation on the right hand side over $ i_2, i_3,i_4, i_5.$  By (\ref{eq:(R2)bianki}), 
$$-(\rr 2\hy 1)_{[i_4 i_5 i_2]\ell_3}{}^{\ell_1}\eta_{[i_3]} = 
(\rr 2\hy 1)_{[i_4 i_5 i_3]\ell_2}{}^{\ell_1}\eta_{[i_2]} =0.$$

For the third equation of (\ref{eq:R2-relation}),
we get 
\begin{align}
\label{eq:R2-2-2}
(\rr 2\hy 2)_{i_2 i_3  \ell_2}{}^{\ell_1} &\times
(\rr 2\hy 2)_{i_4 i_5 \ell_3}{}^{\ell_2} 
\notag\\
&= 
(2(\eta_{i_2} \delta_{i_3}{}^{\ell_1} 
      -\eta_{i_3}\delta_{i_2}{}^{\ell_1})
          \eta_{\ell_2}
 + (g_{i_2 \ell_2} \eta_{i_3} -g_{i_3 \ell_2}\eta_{i_2})
     \xi^{\ell_1})
\notag \\
& \qquad \times 
(2(\eta_{i_4} \delta_{i_5}{}^{\ell_2} 
      -\eta_{i_5}\delta_{i_4}{}^{\ell_2})
          \eta_{\ell_3}
 + (g_{i_4 \ell_3} \eta_{i_5} -g_{i_5 \ell_3}\eta_{i_4})
     \xi^{\ell_2})
\notag \\
&= 2(\eta_{i_2}\delta_{i_3}{}^{\ell_1} 
      -\eta_{i_3}\delta_{i_2}{}^{\ell_1})
   (g_{i_4 \ell_3}\eta_{i_5} -g_{i_5 \ell_3}\eta_{i_4})
 \\
&\qquad 
 +2(\eta_{i_4} \delta_{i_5}{}^{\ell_2} 
      -\eta_{i_5}\delta_{i_4}{}^{\ell_2}) \eta_{\ell_3}
    (g_{i_2 \ell_2} \eta_{i_3} - g_{i_3 \ell_2} \eta_{i_2})
    \xi^{\ell_1}
    \notag \\
&=2(\eta_{i_2}\eta_{i_5} \delta_{i_3}{}^{\ell_1}g_{i_4\ell_3}
-\eta_{i_3}\eta_{i_5}\delta_{i_2}{}^{\ell_1}g_{i_4 \ell_3}
-\eta_{i_2}\eta_{i_4}\delta_{i_3}{}^{\ell_1}g_{i_5 \ell_3}
+\eta_{i_3} \eta_{i_4}\delta_{i_2}{}^{\ell_1}
     g_{i_5 \ell_3})
     \notag \\
& \qquad 
+2(\eta_{i_4}\eta_{i_3}g_{i_5 i_2} 
-\eta_{i_5}\eta_{i_3}g_{i_2 i_4} 
-\eta_{i_4}\eta_{i_2}g_{i_3 i_5}
+\eta_{i_5}\eta_{i_2}g_{i_3 i_4})\xi^{\ell_1}\eta_{\ell_3}.
\notag 
\end{align}
The term $\eta_{i_2}\eta_{i_5} \delta_{i_3}{}^{\ell_1}g_{i_4\ell_3}$ in 
(\ref{eq:R2-2-2}) is symmetric in $i_2, i_5$, so
$$\eta_{[i_2}\eta_{i_5} \delta_{i_3}{}^{\ell_1}g_{i_4]\ell_3} = 
\eta_{i_2}\eta_{i_5}\delta_{[i_3}{}^{\ell_1}g_{i_4]\ell_3} - \eta_{i_5}\eta_{i_2}
\delta_{[i_3}{}^{\ell_1}g_{i_4]\ell_3} = 0.$$
Every other term in the last two lines of (\ref{eq:R2-2-2}) is symmetric in a pair of indices,
so we conclude
\begin{equation*}
((\rr 2\hy 2) \otimes (\rr 2\hy 2))_{[i_2 i_3 i_4 i_5] \ell_3}{}^{\ell_1}  =0.     
\end{equation*}
\end{proof}

By Lemmas \ref{prop:RE-R3} and \ref{lem:5.4},  (\ref{eq:30}) becomes
\begin{equation}\label{eq:aaa}
N_{[i_2 i_3 \cdots i_{4k} i_{4k+1}] }
 {}_{\ell_0} {}^{\ell_1}
= \rho^{4k} 
\underbrace{(
(\rr 2\hy 1)\wedge \cdots \wedge (\rr 2\hy 1)
)}_{\text{2k}}
{}_{i_2 i_3 \cdots i_{4k} i_{4k+1}} 
  {}_{\ell_0} {}^{\ell_1}   
   + l.o.
\end{equation}
For
\begin{equation*}
(\rr 2\hy 1)^\prime_{i_2 i_3 \ell_2}{}^{\ell_1}
=2(\phi_{i_2 \ell_2} \phi_{i_3}{}^{\ell_1}
 +\phi_{i_2 i_3} \phi_{\ell_2}{}^{\ell_1}),
\end{equation*}
we have 
\begin{equation*}
(\rr 2\hy 1)_{[i_2 i_3] \ell_2}{}^{\ell_1} 
=(\rr 2\hy 1)^\prime_{[i_2 i_3] \ell_2}{}^{\ell_1} 
\end{equation*}
Then
 (\ref{eq:aaa}) becomes
\begin{align}  \label{eq:a11} 
N_{[i_2 i_3 \cdots i_{4k} i_{4k+1}]\ell_0}{}^{\ell_1}
&= \rho^{4k}((\rr 2\hy 1)^\prime\wedge
   (\rr 2\hy 1)^\prime 
\cdots \wedge
    (\rr 2\hy 1)^\prime)_{i_2 i_3 \cdots i_{4k} i_{4k+1}}{}_{\ell_0}{}^{\ell_1}.
     \\
      &= \rho^{4k}\underbrace{(
A'\wedge \cdots \wedge A'
)}_{\text{k}}
{}_{i_2 i_3 \cdots i_{4k} i_{4k+1}} 
  {}_{\ell_0} {}^{\ell_1}  
  + l.o.,\notag
\end{align}
where
\begin{align*} 
{A^\prime}_{i_2 i_3 i_4 i_5 \ell_3}{}^{\ell_1}  
&:=
(\rr 2\hy 1)^\prime_{i_2 i_3 \ell_2}{}^{\ell_1}
 \times 
(\rr 2\hy 1)^\prime_{i_4 i_5 \ell_3}{}^{\ell_2} \\
&=4
(\phi_{i_2 \ell_2} \phi_{i_3}{}^{\ell_1}
 +\phi_{i_2 i_3} \phi_{\ell_2}{}^{\ell_1})
 \times
(\phi_{i_4 \ell_3} \phi_{i_5}{}^{\ell_2}
 +\phi_{i_4 i_5} \phi_{\ell_3}{}^{\ell_2})
\notag \\
&=  
4( 
\phi_{i_2 \ell_2}\phi_{i_3}{}^{\ell_1}
  \phi_{i_4 \ell_3}\phi_{i_5}{}^{\ell_2}
+\phi_{i_2 \ell_2}\phi_{i_3}{}^{\ell_1}
 \phi_{i_4 i_5}\phi_{\ell_3}{}^{\ell_2}
\notag \\
&\qquad
+\phi_{i_2 i_3} \phi_{\ell_2}{}^{\ell_1} 
 \phi_{i_4 \ell_3}\phi_{i_5}{}^{\ell_2}
+\phi_{i_2 i_3}\phi_{\ell_2}{}^{\ell_1}
  \phi_{i_4 i_5} \phi_{\ell_3}{}^{\ell_2} 
\notag
) 
\end{align*}
Since $\phi_{i_2 \ell_2}\phi_{i_5}{}^{\ell_2}$ 
is symmetric in $i_2, i_5$,
\begin{equation*}
\phi_{i_2 \ell_2}\phi_{i_3}{}^{\ell_1}
  \phi_{i_4 \ell_3}\phi_{i_5}{}^{\ell_2}
\sim
0  \quad \modd  (i_2 i_3 i_4 i_5). 
\end{equation*}
 Since exchanging the indices $i_2, i_3$ with $i_4, i_5$
does not change the sign of the permutation, 
\begin{equation*}
\phi_{i_2 \ell_2}\phi_{i_3}{}^{\ell_1}
  \phi_{i_4 i_5}\phi_{\ell_3}{}^{\ell_2}
\sim
\phi_{i_4 \ell_2}\phi_{i_5}{}^{\ell_1}
 \phi_{i_2 i_3}\phi_{\ell_3}{}^{\ell_2}
\quad \modd (i_2 i_3 i_4 i_5). 
\end{equation*}
Therefore,
\begin{align*}\label{eq:demod} 
&{A^\prime}_{i_2 i_3 i_4 i_5 \ell_3}{}^{\ell_1}
\\
&\qquad \sim
4(
\phi_{i_4 \ell_2} 
    \phi_{\ell_3}{}^{\ell_2}\phi_{i_5}{}^{\ell_1} 
+
\phi_{i_4 \ell_3} 
    \phi_{i_5}{}^{\ell_2}\phi_{\ell_2}{}^{\ell_1} 
+
\phi_{i_4 i_5} 
    \phi_{\ell_2}{}^{\ell_1}\phi_{\ell_3}{}^{\ell_2} )
    \phi_{i_2 i_3}
\quad \modd  (i_2 i_3 i_4 i_5), \notag
\end{align*}
which we will use to compute $(A'\wedge\ldots\wedge A'){}_{i_2 i_3\cdots i_{4k} i_{4k+1} \ell}{}^{\ell}$ in (\ref{eq:a11}).
Since the wedge product involves skew-symmetrization, or modding out in our notation, 
we can relabel $A'$  as
\begin{equation}\label{??}
(A')_{i_2 i_3 i_4 i_5 \ell_3}{}^{\ell_1} 
:= 4\phi_{i_2 i_3}(B^*)_{i_4 i_5 \ell_3}{}^{\ell_1}
:=4\phi_{i_2 i_3}((\rb 1)_{i_4 i_5 \ell_3}{}^{\ell_1}+(\rb 2)_{i_4 i_5 \ell_3}{}^{\ell_1}
+(\rb 3)_{i_4 i_5 \ell_3}{}^{\ell_1} ),
\end{equation}
where 
\begin{align}\label{b*}
(B^*)_{i_4 i_5 \ell_3}{}^{\ell_1} &= (\rb 1)_{i_4 i_5 \ell_3}{}^{\ell_1}+
(\rb 2)_{i_4 i_5 \ell_3}{}^{\ell_1}+ (\rb 3)_{i_4 i_5 \ell_3}{}^{\ell_1},\\
(\rb 1)_{i_4 i_5 \ell_3}{}^{\ell_1}
&=-(\phi^{(2)})_{i_4 \ell_3} \phi_{i_5}{}^{\ell_1},
\ (\rb 2)_{i_4 i_5 \ell_3}{}^{\ell_1}
=\phi_{i_4 \ell_3}(\phi^{(2)})_{i_5}{}^{\ell_1},
\ (\rb 3)_{i_4 i_5 \ell_3}{}^{\ell_1} 
=\phi_{i_4 i_5} (\phi^{(2)})^{\ell_1}{}_{\ell_3}. \notag
\end{align}
Here we introduce the notation
\begin{align*}
(\phi^{(2)})_k{}^h &:=     \phi_k{}^{j_1} \phi_{j_1}{}^h   = -\delta_{k}{}^h + \eta_{k} \xi^{h},\ 
(\phi^{(3)})_k{}^h :=     \phi_k{}^{j_1} \phi_{j_1}{}^{j_2} \phi_{j_2}{}^h  = -\phi_{k}{}^{h}, \\
(\phi^{(4)})_k{}^h &:=  \phi_k{}^{j_1} \phi_{j_1}{}^{j_2}\phi_{j_2}{}^{j_3} \phi_{j_3}{}^h 
    = \delta_{k}{}^{h} - \eta_{k} \xi^{h}    =-(\phi^{(2)})_k{}^{h}.    
\end{align*}
For
\begin{equation}\label{!!!}
(A')^{(2)}_{i_2 \cdots i_9 \ell_4}{}^{\ell_1}
:=(A')_{i_2 i_3 i_4 i_5 \ell_3}{}^{\ell_1} \times
(A')_{i_6 i_7 i_8 i_9 \ell_4}{}^{\ell_3} 
 =4^2 \phi_{i_2 i_3}\phi_{i_6 i_7} (B^*)_{i_4 i_5 \ell_3}{}^{\ell_1}\times
 (B^*)_{i_8 i_9 \ell_4}{}^{\ell_3}
\end{equation}
we have:

\begin{lem} \label{lem:B^2}
\begin{equation*}
(A')^{(2)}_{i_2 \cdots i_9 \ell_4}{}^{\ell_1}
\sim
-(4^2) \phi_{i_2 i_3} \phi_{i_4 i_5} \phi_{i_6 i_7}
B^*_{i_8 i_9 \ell_4}{}^{\ell_1}
\quad \modd (i_2 \ldots i_9). 
\end{equation*}
\end{lem}
\begin{proof} We compute
\begin{align*}
(\rb 1)_{i_4 i_5\ell_3}{}^{\ell_1} 
\times 
(\rb 1)_{i_8 i_9\ell_4}{}^{\ell_3} &=
(\phi^{(2})_{i_4 \ell_3} \phi_{i_9}{}^{\ell_3} 
(\phi^{(2)})_{i_8 \ell_4} \phi_{i_5}{}^{\ell_1}
=
-(\phi^{(3)})_{i_4 i_9} (\phi^{(2)})_{i_8 \ell_4} \phi_{i_5}{}^{\ell_1}
\notag \\
& \sim 
\phi^{(3)}_{i_4 i_5} \phi^{(2)}_{i_8 \ell_4} \phi_{i_9}{}^{\ell_1} 
\sim -\phi_{i_4 i_5}\phi^{(2)}_{i_8 \ell_4} \phi_{i_9}{}^{\ell_1} 
\quad \modd  (i_4 i_5 i_8 i_9),\\
(\rb 1)_{i_4 i_5\ell_3}{}^{\ell_1} 
\times 
(\rb 2)_{i_8 i_9\ell_4}{}^{\ell_2} 
&=
-(\phi^{(2)})_{i_4 \ell_3}\phi_{i_5}{}^{\ell_1}
  \phi_{i_8 \ell_4}(\phi^{(2)})_{i_9}{}^{\ell_3}
=-(\phi^{(4)})_{i_4 i_9} \phi_{i_5}{}^{\ell_1}  \phi_{i_8 \ell_4}\notag \\
& \sim 0 \quad \modd  (i_4 i_5 i_8 i_9),\\
(\rb 1)_{i_4 i_5\ell_3}{}^{\ell_1} \times (\rb 3)_{i_8 i_9\ell_4}{}^{\ell_2} 
& = -(\phi^{(2)})_{i_4 \ell_3} \phi_{i_5}{}^{\ell_1}
(\phi^{(2)})^{\ell_3}{}_{\ell_4}\phi_{i_8 i_9}
=-(\phi^{(4)}_{i_4 \ell_4}) \phi_{i_5}{}^{\ell_1}\phi_{i_8 i_9}\\
&\sim -(\phi^{(4)})_{i_8 \ell_4} \phi_{i_9}{}^{\ell_1}\phi_{i_4 i_5}
 \sim (\phi^{(2)})_{i_8 \ell_4}\phi_{i_9}{}^{\ell_1} \phi_{i_4 i_5}
\quad \modd (i_4 i_5 i_8 i_9).
\end{align*}
(For the second equation, note that 
$\phi^{(4)}_{i_4 i_9}    = g_{i_4 i_9} -\eta_{i_4}\eta_{i_9}$
is symmetric in $i_4, i_9$.)

Thus, we have
\begin{equation}
\label{eq:B1-BE}\sum_{E=1}^3 (\rb 1)_{i_4 i_5 \ell_3}{}^{\ell_1}  
\times (\rb {\text E})_{i_8 i_9 \ell_4}{}^{\ell_3}  
\sim 0 \quad \modd (i_4 i_5 i_8 i_9).
\end{equation}

Again using the symmetry of $(\phi^{(2)})_{i_ai_b}$, we have 
\begin{align*}
(\rb 2)_{i_4 i_5\ell_3}{}^{\ell_1} \times (\rb 1)_{i_8 i_9\ell_4}{}^{\ell_2} 
&= -\phi_{i_4 \ell_3}(\phi^{(2)})_{i_5}{}^{\ell_1}
     (\phi^{(2)})_{i_8 \ell_4}
      \phi_{i_9}{}^{\ell_3}
= (\phi^{(2)})_{i_8 \ell_4} 
   (\phi^{(2)})_{i_4 i_9} 
   (\phi^{(2)})_{i_5}{}^{\ell_1} \\
&\sim 0 \quad \modd  (i_4 i_5 i_8 i_9),\\
(\rb 2)_{i_4 i_5\ell_3}{}^{\ell_1} \times 
(\rb 2)_{i_8 i_9\ell_4}{}^{\ell_2} &\\  
= 
\phi_{i_4 \ell_3}(\phi^{(2)})_{i_5}{}^{\ell_1}&
 \phi_{i_8 \ell_4} (\phi^{(2)})_{i_9}{}^{\ell_3}
=-(\phi^{(3)})_{i_9 i_4}    (\phi^{(2)})_{i_5}{}^{\ell_1}  \phi_{i_8 \ell_4}
=\phi_{i_9 i_4}(\phi^{(2)})_{i_5}{}^{\ell_1}  \phi_{i_8 \ell_4}
\notag \\
 \sim
-\phi_{i_5 i_4}(\phi^{(2)})_{i_9}&{}^{\ell_1}
  \phi_{i_8 \ell_4} \sim
\phi_{i_4 i_5}(\phi^{(2)})_{i_9}{}^{\ell_1}
  \phi_{i_8 \ell_4}
\quad \modd  (i_4 i_5 i_8 i_9),\\
(\rb 2)_{i_4 i_5\ell_3}{}^{\ell_1} \times (\rb 3)_{i_8 i_9\ell_4}{}^{\ell_2}   
&= 
\phi_{i_4 \ell_3}(\phi^{(2)})_{i_5}{}^{\ell_1}
   \phi_{i_8 i_9}(\phi^{(2)})^{\ell_3}{}_{\ell_4}
=\phi_{i_8 i_9} \phi_{i_4 \ell_3}(\phi^{(2)})^{\ell_3}{}_{\ell_4}   (\phi^{(2)})_{i_5}{}^{\ell_1}\notag \\
&=
\phi_{i_8 i_9} (\phi^{(3)})_{i_4 \ell_4}
  (\phi^{(2)})_{i_5}{}^{\ell_1}
\notag \\    
&\sim 
\phi_{i_4 i_5} (\phi^{(3)})_{i_8 \ell_4}      (\phi^{(2)})_{i_9}{}^{\ell_1}    
\sim 
-\phi_{i_4 i_5} \phi_{i_8 \ell_4} 
     (\phi^{(2)})_{i_9}{}^{\ell_1}
\quad \modd  (i_4 i_5 i_8 i_9).
\end{align*}
Therefore,
\begin{equation}
\label{eq:B2-BE}
\sum_{E=1}^3 
(\rb 2)_{i_4 i_5 \ell_3}{}^{\ell_1}  \times (\rb \text{E})_{i_8 i_9 \ell_4}{}^{\ell_3}  
\sim 0 \quad \modd (i_4 i_5 i_8 i_9).
\end{equation}

 Finally, we compute 
\begin{align}
\label{eq:B3-B1}
(\rb 3)_{i_4 i_5\ell_3}{}^{\ell_1} \times (\rb 1)_{i_8 i_9\ell_4}{}^{\ell_3} 
&=
\phi_{i_4 i_5} (\phi^{(2)})^{\ell_1}{}_{\ell_3}
(-(\phi^{(2)})_{i_8 \ell_4} )\phi_{i_9}{}^{\ell_3}
=\phi_{i_4 i_5} (\phi^{(3)})^{\ell_1}{}_{i_9}  (\phi^{(2)})_{i_8 \ell_4}
\notag \\
&=
\phi_{i_4 i_5}(-\phi^{\ell_1}{}_{i_9} )
   (\phi^{(2)})_{i_8 \ell_4} =\phi_{i_4 i_5}\phi_{i_9}{}^{\ell_1}   (\phi^{(2)})_{i_8 \ell_4}\notag\\
& \sim   -\phi_{i_4 i_5}\phi_{i_8}{}^{\ell_1}   (\phi^{(2)})_{i_9 \ell_4}
\quad \modd  (i_4 i_5 i_8 i_9),\notag\\
(\rb 3)_{i_4 i_5\ell_3}{}^{\ell_1} \times (\rb 2)_{i_8 i_9\ell_4}{}^{\ell_3} 
&=  \phi_{i_4 i_5} (\phi^{(2)})^{\ell_1}{}_{\ell_3}
   \phi_{i_8 \ell_4} (\phi^{(2)})_{i_9}{}^{\ell_3}
= \phi_{i_4 i_5} \phi_{i_8 \ell_4}    (\phi^{(4)})^{\ell_1}{}_{i_9}     \\
&=-\phi_{i_4 i_5} \phi_{i_8 \ell_4}    (\phi^{(2)})_{i_9}{}^{\ell_1},\notag\\
(\rb 3)_{i_4 i_5\ell_3}{}^{\ell_1} \times (\rb 3)_{i_8 i_9\ell_4}{}^{\ell_3} 
&=\phi_{i_4 i_5} \phi_{i_8 i_9}(\phi^{(2)})^{\ell_1}{}_{\ell_3}
   (\phi^{(2)})^{\ell_3}{}_{\ell_4} = \phi_{i_4 i_5} \phi_{i_8 i_9}
   (\phi^{(4)})^{\ell_1}{}_{\ell_4}\notag\\
   &= -\phi_{i_4 i_5} \phi_{i_8 i_9}(\phi^{(2)})^{\ell_1}{}_{\ell_4}.\notag
\end{align}

By (\ref{eq:B1-BE}), (\ref{eq:B2-BE}),
(\ref{eq:B3-B1}), we have 
\begin{align*}
(A')^{(2)}_{i_2 \cdots i_9 \ell_4}{}^{\ell_1}
& \sim \phi_{i_2 i_3} \phi_{i_6 i_7}
\sum_{E=1}^3 
  (\rb 3)_{i_4 i_5 \ell_3}{}^{\ell_1}  \times  (\rb \text{E})_{i_8 i_9 \ell_4}{}^{\ell_3}
  \quad \modd  (i_2\ldots i_9)\\
 &\sim -(4^2) \phi_{i_2 i_3} \phi_{i_4 i_5} \phi_{i_6 i_7}
B^*_{i_8 i_9 \ell_4}{}^{\ell_1} \quad \modd  (i_2\ldots i_9),     
\end{align*}
which proves Lemma \ref{lem:B^2}. 
 \end{proof}

\begin{cor}\label{cor:quick5} For $(A')_{i_2 i_3 i_4 i_5 \ell_3}{}^{\ell_1} $ in (\ref{??}) and 
$(A')^{(2)}_{i_2 \cdots i_9 \ell_4}{}^{\ell_1}$ in (\ref{!!!}), set
$$(A')^{(r)}_{i_2\ldots i_{4r+1}\ell_{r+2}}{}^{\ell_1} = (A')^{(r-1)}_{i_2\ldots i_{4r-3}\ell_{r+1}}{}^{\ell_1}(A')_{i_{4r-2}i_{4r-1}i_{4r}i_{4r+1}\ell_{r+2}}{}^{\ell_{k+1}},$$
$r= 2,\ldots, k.$
Then 
(i)
$$(A')^{(r)}_{i_2\ldots i_{4r+1}\ell_{r+2}}{}^{\ell_1} \sim (-1)^{r-1} 4^r \phi_{i_2i_3}\phi_{i_4i_5}
\ldots \phi_{i_{4r-2} i_{4r-1} } B^*_{i_{4r}i_{4r+1} \ell_{r+2}} {}^{\ell_1}
\quad \modd(i_2\ldots i_{4r+1})$$;

(ii) In the notation of (\ref{prop:N}),
$$\tr(N)_{[i_2 \cdots i_{4k+1}]} =  N_{[i_2 \cdots i_{4k+1}]}{}_{\ell}{}^{\ell} 
= (-1)^{k-1} 4^k (4k+2) \rho^{4k}(\phi\wedge\ldots\wedge \phi)_{i_2\ldots i_{4k+1}} + l.o.$$
\end{cor}

\begin{proof}
(i) By Lemma \ref{lem:B^2}, this follows by induction on $r$ up to $r=k$, where dim$(M)= 4k+1.$

(ii) By (\ref{prop:N}), (\ref{eq:regroup}), (\ref{eq:a11}), (\ref{??}), 
\begin{align}\label{trn}
    \tr(N)_{[i_2 \cdots i_{4k+1}]} &= \rho^{4k}\tr (A')^{(k)}_{[i_2\ldots i_{4k+1}]}
    = \rho^{4k}(A')^{(k)}_{[i_2\ldots i_{4k+1}]\ell} {}^{\ell}\\
    &\sim (-1)^{k-1} 4^k\rho^{4k} \phi_{[i_2i_3}\phi_{i_4i_5}
\ldots \phi_{i_{4k-2} i_{4k-1} } B^*_{i_{4k}i_{4k+1}] \ell} {}^{\ell}\quad \modd (i_2\ldots i_{4k+1}).\notag
\end{align}
By (\ref{b*}), 
\begin{align*}B^*_{ab \ell} {}^{\ell} &= -(\phi^{(2)})_{a \ell} \phi_{b}{}^{\ell}+
\phi_{a \ell}(\phi^{(2)})_{b}{}^{\ell}
+\phi_{ab} (\phi^{(2)})^{\ell}{}_{\ell} = \phi_{ab}+\phi_{ab}+ \phi_{ab} (\delta^\ell{}_\ell - \eta_\ell\xi^\ell)\\
&= 2\phi_{ab}+\phi_{ab} ((4k+1) - 1) = (4k+2)\phi_{ab},
\end{align*}
and plugging this into (\ref{trn}) finishes the proof.
\end{proof}

This proves (\ref{prop:N}),
and so finishes the proof of Prop.~\ref{prop:1}.

\section{Proof of Lemma \ref{lem:5}}\label{app:S6}

\begin{customlem}{6.1}  For $(S^{4k+1}, h_\rho)$, in the notation of (\ref{123}) we have
\begin{align*}
(i) &\qquad 
(({\rm R}1)\otimes ({\rm R} {\rm E}))_{[i_2 i_3 i_4 i_5] \ell_3}{}^{\ell_1}   
=0;
\\
(ii) &\qquad 
(({\rm R}{\rm E})\otimes ({\rm R}1))_{[i_2 i_3 i_4 i_5] \ell_3}{}^{\ell_1}   
=0;
\end{align*}
for ${\rm E}=1,2,3.$
\end{customlem}
\begin{proof} (i)
$(S^{4k+1}, g_{std})$ has positive constant sectional curvature $g(R(X,Y)Y,X) = 1$, for $|X|=|Y|=1, X\perp Y$, and curvature tensor
\begin{equation}
\label{constant-curvature form}
  R_{kji}{}^h 
   = \frac{c+3}{4}(g_{ki}\delta_j{}^h - g_{ji}\delta_k{}^h),
\end{equation}
with $c=1$ \cite[I.V.2]{K-N}. 

Noting that  
$(\rr 1)_{kji}{}^{h} = R_{kji}{}^{h}, $
we have  
\begin{align*}
(\rr 1)_{i_2 i_3 \ell_2}{}^{\ell_1} \times
( \re)_{i_4 i_5 \ell_3}{}^{\ell_2}
&= \frac{c+3}{4}
(g_{i_2 \ell_2}\delta_{i_3}{}^{\ell_1}
 -g_{i_3 \ell_2} \delta_{i_2}{}^{\ell_1})
 \times (\re)_{i_4 i_5 \ell_3}{}^{\ell_2}
  \\
 & 
 = \frac{c+3}{4}
 (\delta_{i_3}{}^{\ell_1}(\re)_{i_4 i_5 \ell_3 i_2}
   -\delta_{i_2}{}^{\ell_1}(\re)_{i_4 i_5 \ell_3 i_3}),
   \notag 
\end{align*}
where
$(\re)_{kjih} = g_{\ell h} (\re)_{kji}{}^{\ell}.$

Note that 
$(\re)_{kjih}=- (\re)_{kjhi}:$ for $\bar R_{kjih} = -\bar R_{kjhi}$ for any $\rho$, and $(\rr 1), (\rr 2), (\rr 3)$ come with different powers of $\rho$, so each $(\re)$ is skew-symmetric in $i,h.$
Thus, 
\begin{equation*}
((\rr 1) \otimes (\RE))_{[i_2 i_3 i_4 i_5] \ell_3}{}^{\ell_1} 
= \frac{c+3}{4}
 (-\delta_{[i_3}{}^{\ell_1}(\RE)_{i_4 i_5 i_2] \ell_3 }
   -\delta_{[i_2}{}^{\ell_1}(\RE)_{i_4 i_5 i_3] \ell_3})
\end{equation*}
where we skew-symmetrize over the $i$ indices. 
Since $(\RE)_{[kji]}{}^{h}=0$ by the Bianchi identity,
\begin{equation*}
((\rr 1) \otimes (\RE))_{[i_2 i_3 i_4 i_5] \ell_3}{}^{\ell_1}  =0.   
\end{equation*}

(ii) This is similar:
\begin{align*}
 (\RE)_{i_2 i_3 \ell_2}{}^{\ell_1}  
 \times (\rr 1)_{i_4 i_5 \ell_3}{}^{\ell_2}
 &=(\RE)_{i_2 i_3 \ell_2}{}^{\ell_1}
 \times
 \frac{c+3}{4}
(g_{i_4 \ell_3}\delta_{i_5}{}^{\ell_2}
 -g_{i_5 \ell_3} \delta_{i_4}{}^{\ell_2})\\
& =
 \frac{c+3}{4} 
 ((\RE)_{i_2 i_3 i_5}{}^{\ell_1} g_{i_4 \ell_3}
 -
(\RE)_{i_2 i_3 i_4}{}^{\ell_1} g_{i_5 \ell_3} ),
\notag 
\end{align*}
so again,
\begin{equation*}
((\RE) \otimes (\rr 1))_{[i_2 i_3 i_4 i_5] \ell_3}{}^{\ell_1}  =0.   
\end{equation*}
\end{proof}

\bibliographystyle{amsplain}
\bibliography{arXivIV}

\end{document}